\documentclass[11pt]{article}
\usepackage[T1]{fontenc}
\usepackage{tikz}
\usepackage{lmodern}
\usepackage{amsmath}
\usepackage{amsthm}
\usepackage{amssymb}
\usepackage{graphicx}
\usepackage{multirow}
\usepackage{multicol}
\usepackage{booktabs}
\usepackage[T1]{fontenc}
\usepackage{microtype}
\usepackage{floatrow}
\usepackage[english]{babel}
\usepackage[hmarginratio=1:1,top=32mm,columnsep=20pt]{geometry}
\usepackage[hang, small,labelfont=bf,up,textfont=it,up]{caption}
\usepackage{indentfirst}
\usepackage{diagbox}
\usepackage{subfigure}
\usepackage{abstract}

\usepackage{titlesec}
\usepackage{titling}
\usepackage{hyperref}
\usepackage[ruled]{algorithm2e}
\usepackage{tkz-euclide,pgfplots}
\newtheorem{lemma}{Lemma}
\newtheorem{theorem}{Theorem}

\newtheorem{definition}{Definition}
\newtheorem{remark}{Remark}

\pretitle{\begin{center}\huge\bfseries}
\posttitle{\end{center}}

\title{Robust covariance estimation for distributed principal component analysis\thanks{This work was supported by grants from the NSF of China (Grant No.11731012), Ten Thousands Talents Plan of Zhejiang Province (Grant No. 2018R52042) and the Fundamental Research Funds for the Central Universities.}  }

\author{Kangqiang Li\thanks{Corresponding author E-mail address: 11935023@zju.edu.cn (Kangqiang Li)} \qquad Han Bao\thanks{E-mail address: 21935024@zju.edu.cn (Han Bao)}\qquad Lixin Zhang\thanks{E-mail address: stazlx@zju.edu.cn (Lixin Zhang)}\\[1ex] % Your name
\normalsize School of Mathematical Sciences, Zhejiang University,  Hangzhou,  Zhejiang 310027,  China \\ % Your institution
}

\begin{document}

% Print the title
\maketitle
\begin{abstract}
Fan et al. \newblock{[\em{Annals of Statistics}} \textbf{47}(6) (2019) 3009--3031] constructed a distributed principal component analysis (PCA) algorithm to reduce the communication cost between multiple servers significantly. However, their algorithm's guarantee is only for sub-Gaussian data. Spurred by this deficiency, this paper enhances the effectiveness of their distributed PCA algorithm by utilizing robust covariance matrix estimators of Minsker \newblock{[\em{Annals of Statistics}} \textbf{46}(6A) (2018) 2871--2903] and Ke et al. \newblock{[\em{Statistical Science}} \textbf{34}(3) (2019) 454--471] to tame heavy-tailed data. The theoretical results demonstrate that when the sampling distribution is symmetric innovation with the bounded fourth moment or asymmetric with the finite $6$-th moment, the statistical error rate of the final estimator produced by the robust algorithm is similar to that of sub-Gaussian tails. Extensive numerical trials support the theoretical analysis and indicate that our algorithm is robust to heavy-tailed data and outliers.
\end{abstract}

%----------------------------------------------------------------------------------------
%	ARTICLE CONTENTS
%----------------------------------------------------------------------------------------

\section{\textbf{Introduction}}\label{sectIntr}
Principal component analysis (PCA) is a well-known dimensionality reduction technique initiated by Pearson (1901)\cite{pearson1901liii} that extracts from the original data the main information and discards the irrelevant information. In early research, Davis (1977)\cite{doi:10.1111/j.1467-842X.1977.tb01088.x} extended Anderson's (1963)\cite{anderson1963} asymptotic results about PCA under Gaussian assumption to the case of non-Gaussian distributions with only bounded fourth moment. Davis derived that when the fourth-order moment of samples exists, the eigenvalues and eigenvectors of the sample covariance matrix enjoy a similar statistical error rate for the case of Gaussian distribution. Although the sample covariance matrix has been applied in numerous fields, in the theoretical properties, it relies heavily on Gaussian distribution (Wang and Fan (2017)\cite{wang2017asymptotics}) and presents poor performance if the Gaussian assumption is neglected (Bickel and Levina (2008)\cite{bickel2008covariance}, Fan, Fan and Lv (2008)\cite{fan2008high}). Additionally, the conventional PCA does not consider the presence of outliers and heavy-tailed data. Therefore, constructing tail-robust covariance matrix estimators with heavy-tailed data has been a hot research topic over the past few decades. For example, regarding the elliptically symmetric distribution, which is a popular way to simulate heavy-tailed data, Fan, Liu and Wang (2018)\cite{Large covariance} designed a procedure to estimate the spiked covariance matrix of the elliptical factor model when the random samples have the fourth moment bounded. They derived the error bound of the robust estimator under max norm and the statistical error rate of the corresponding eigenvalues and eigenvectors. Furthermore, Fan, Wang and Zhong (2019)\cite{approximate factor} employed the approximate factor model to recover the covariance matrix of the heavy-tailed data. Fan, Wang and Zhu (2021)\cite{A shrinkage principle} constructed a $\ell_{4}$-norm truncated sample covariance matrix when the bounded fourth moment of heavy-tailed distributions exists, and they demonstrated that the robust covariance estimator has the statistical error rate of $O_{P} (\sqrt{d\log d/n})$ under the spectral norm. Under only the $(2+\epsilon)$-th moment assumption, Avella-Medina et al. (2018)\cite{precision matrices} proposed employing Huber (1964)\cite{location parameter}'s loss and adaptive thresholding function to construct robust estimators for recovering weakly sparse covariance and precision matrix. In this study, the phase transition for the convergence rate of the covariance estimation was first observed. From the random matrix perspective, Minsker (2018)\cite{Minsker_2018} and  Ke et al. (2019)\cite{User-friendly}, respectively, proposed different robust covariance matrix estimators such that the computing amount of estimating the population covariance matrix can be remarkably reduced and derived exponential-type concentration inequality under the spectral norm. Minsker generalized Catoni's (2012)\cite{catoni2012} $M$-estimator from a robust mean estimation of the univariate random variables to a multivariate version of the self-adjoint random matrices. Ke employed Huber's (1964)\cite{{location parameter}} loss function and its derivative to construct a Huber-type $M$-estimator and a truncated covariance matrix estimator for heavy-tailed distributions with the bounded fourth moment. Finally, it should be noted that a lot of literature suggests several robust covariance estimators such as Catoni (2016)\cite{PAC-Bayesian}, Minsker and Wei (2020)\cite{Robust modifications} and Mendelson and Zhivotovskiy (2018)\cite{Robust covariance}.\par

With the emergence of massive datasets, several distributed algorithms based on various statistical models have been recently suggested including Tian and Gu (2017)\cite{Communication4} for distributed sparse linear discriminant analysis, Lee et al. (2017)\cite{Communication1}, Jordan et al. (2018)\cite{Communication3} and Fan, Guo and Wang (2019)\cite{Communication-efficient} for distributed regression model, and El Karoui and d'Aspremont (2010)\cite{el2010second}, Schizas and Aduroja (2015)\cite{Schizas}, Chen et al (2016)\cite{chen2016integrating}, Fan et al. (2019)\cite{fan2019distributed} and Chen et al (2021)\cite{Chen} for distributed PCA. Fan et al. (2019)\cite{fan2019distributed} constructed a distributed PCA algorithm via utilizing the sample covariance matrix. Specifically, they computed the leading eigenvectors of the sample covariance matrix on each server and forwarded them to the data center. Their theoretical analysis and numerical simulation assumed that the samples obey the sub-Gaussian distribution, which in practice is very restrictive since real-world data usually do not exactly satisfy this assumption. To address this problem, in this paper, in the case of finite $(2+\epsilon)$-th moment, we first focus on the tail-robust covariance matrix estimators proposed by Minsker (2018)\cite{Minsker_2018} and Ke et al. (2019)\cite{User-friendly}, and then we enhance the robustness of the distributed PCA algorithm by utilizing these user-friendly covariance estimators to ensure the algorithm's reliability in the presence of heavy-tailed data. Even when some servers have many large outliers, the robust algorithm can still be useful. Besides, the theoretical results show that the distributed estimator achieves a similar statistical error rate to that of sub-Gaussian distributed data in Fan et al. (2019)\cite{fan2019distributed}, when the sampling distribution is symmetric innovation and has only finite fourth moment. For the generic heavy-tailed distribution, the proposed algorithm applies to the finite $6$-th moment case. We only require samples with symmetric innovation satisfying the bounded $4$-th moment condition, which makes it easier for the occurrence of data far away from the mean. Our simulation study supports the corresponding theoretical analysis.\par

The remainder of this paper is organized as follows: Section \ref{SecDef}  introduces some definitions and notation, while Section \ref{SecMeth} introduces the proposed robust distributed PCA algorithm. Section \ref{SecMain} revisits two user-friendly robust covariance matrices and generalizes their application conditions, and we embed them in the suggested algorithm. Extensive simulation results are presented in Section \ref{SecSimulation} supporting the theoretical analysis. A brief discussion is provided in Section \ref{SecDiscussion}. Finally, all the proofs are presented in Appendix A.

%------------------------------------------------
\section{\textbf{Notation and definitions}}\label{SecDef}
Given $A \in \mathbb{R}^{d \times d}, $ let $A^{T} \in \mathbb{R}^{d \times d}$ be the transpose of $A$ and $\operatorname{tr}(A)$ be the trace of $A$. The matrix norms are defined as follows: the spectral and Frobenius norms are denoted by $\|A\|_{2}:=\sqrt{\lambda_{\max }\left(A^{T} A\right)}$ and $\|A\|_{\mathrm{F}}=\sqrt{\operatorname{tr}\left(A^{T} A\right)}$, respectively, where $\lambda_{\max}(A^{T}A)$ is the largest eigenvalue of $A^{T}A$. If $A$ is symmetric, the max norm is defined as $\|A\|_{\max}=\max_{1\leq i,j \leq d}|A_{(i,j)}|$. For two vectors $X, Y \in \mathbb{R}^{d}$, $\|X\|_{2}:=\sqrt{X^{T}X}$ and $\left\langle X, Y\right\rangle:=X^{T}Y$. For any positive integer $n$, the set $\{1,2,\ldots,n\}$ is denoted by $[n]$. Given two sequences $\{a_{n}\}_{n=1}^{\infty}$ and $\{b_{n}\}_{n=1}^{\infty}$, if there are some positive constants $C, C_{1}$ such that $Cb_{n}\leq a_{n}\leq C_{1}b_{n}$ holds for $n\geq 1$, we denote it by $a_{n}=O(b_{n})$, and if only $a_{n}\leq C_{1}b_{n}$ holds, we denote it by $a_{n}\lesssim b_{n}$. For a random variable $X \in \mathbb{R}^{d},$ we denote $X_{(k)}$ as the $k$-th entry of $X$, and if $d=1$, we define $\psi_{1}$-norm of $X$ as $\|X\|_{\psi_{1}}=\sup _{k \geq 1}\sqrt[k]{\mathbb{E}|X|^{k}} / k$. $T_{d}(\nu, \mu, \Sigma)$ is the  $d$-variate $t$ distribution with degrees of freedom $\nu, $ mean vector $\mu$ and positive definite matrix $\Sigma$. In this paper, we suppose that the random sample  $X \in \mathbb{R}^{d}$ has the population covariance matrix $\Sigma$ and $\Delta:=\lambda_{K}-\lambda_{K+1}$ is positive such that the top $K$ eigenspace is uniquely determined. $r:=\operatorname{tr}(\Sigma) / \lambda_{1}$ is the effective rank of $\Sigma$.\par

The following definition is presented in Minsker (2018)\cite{Minsker_2018} and Ke et al. (2019)\cite{User-friendly}.
\begin{definition}\label{definition1}  Let $f$ be a real-valued function defined on $\mathbb{R}, $ and $A \in \mathbb{R}^{d \times d}$ be a symmetric matrix with the eigenvalue decomposition $A=V \Lambda V^{T}$ such that $\lambda_{i}(A) \in \mathbb{R}$, $i=1,  \ldots,  d .$ We define $f(A)$ as $f(A)=V f(\Lambda) V^{T}, $ where $f(\Lambda)=\operatorname{diag}(f\left(\lambda_{1}\right), f\left(\lambda_{2}\right), \ldots, f\left(\lambda_{d}\right))$.
\end{definition}

The following definition is defined in Fan et al. (2019)\cite{fan2019distributed}.
\begin{definition} \label{definition3} For a $d$-dimensional random vector  $X$ with covariance  $\Sigma=V \Lambda V^{T}$, if for all $j \in[d]$,  $Z \stackrel{d}{=}\left(I_{d}-2 e_{j} e_{j}^{T}\right)Z$ where  $Z=\Lambda^{-\frac{1}{2}} V^{T} X$, $X$ is named as symmetric innovation.
\end{definition}
The symmetric innovation embodies various distributions such as multivariate Gaussian distribution, multivariate $t$ distribution with mean zero, and their generalization: symmetric elliptical distribution.

\section{\textbf{Robust distributed PCA algorithm}}\label{SecMeth}
\begin{algorithm}[h]
\caption{Robust Distributed PCA Algorithm}
\label{Al}
\KwIn{Data $\left\{X_{i}^{(\ell)}\right\}_{i=1,\ell=1}^{n,m}$ with $\left\|\mathbb{E}|X_{i}^{(\ell)}{X_{i}^{(\ell)T}}|^2\right\|_{2}<\infty$, the tuning parameter $\theta_{(\ell)}$ (or $\tau_{(\ell)}$).}
\For{$\ell=1,2,3,\ldots,m$}{
The $\ell$-th server calculates locally the robust covariance matrix estimator $\widehat{\Sigma}_{n}^{(\ell)}(2, \theta_{(\ell)})$ and its top $K$ eigenvectors
$\widehat{V}_{K}^{(\ell)}=\left(\widehat{v}_{1}^{(\ell)},  \ldots,  \widehat{v}_{K}^{(\ell)}\right) \in \mathbb{R}^{d \times K}$, and transmits the eigenvectors to the central server.
}
\If{Receive $\left\{\widehat{V}_{K}^{(\ell)}\right\}_{\ell=1}^{m}$ from $m$ servers}{
The central server calculates $\widetilde{\Sigma}=\frac{1}{m} \sum_{\ell=1}^{m} \widehat{V}_{K}^{(\ell)} \widehat{V}_{K}^{(\ell)^{T}}$ and its top $K$ eigenvectors $\left\{\widetilde{v}_{i}\right\}_{i=1}^{K}$.}
\KwOut{$\widetilde{V}_{K}=\left(\widetilde{v}_{1}, \ldots, \widetilde{v}_{K}\right) \in \mathbb{R}^{d \times K}$.
}
\end{algorithm}
Let $N=m \cdot n$  be independent samples $\left\{X_{i}^{(\ell)}\right\}_{i=1,\ell=1}^{n,m}$ that are dispersed on $m$ servers with $\mathbb{E} X_{i}^{(\ell)}=0_{d \times 1}, \mathbb{E} X_{i}^{(\ell)}{X_{i}^{(\ell)T}}=\Sigma^{(\ell)}$ and $\left\|\mathbb{E}|X_{i}^{(\ell)}{X_{i}^{(\ell)T}}|^{2}\right\|_{2}<\infty$ for all $i \in [n]$. $\{\Sigma^{(\ell)}\}_{\ell=1}^{m}$ have the same top $K$ eigenvectors $V_{K}$ and each server $\ell$ possesses $n$ samples. In the heavy-tailed setting, we plug the user-friendly robust covariance matrix estimators of Minsker (2018)\cite{Minsker_2018} and Ke et al. (2019)\cite{User-friendly} in the proposed robust distributed PCA algorithm.\par

Similar to Section 4 in Fan et al. (2019)\cite{fan2019distributed}, to measure the statistical error $\rho\left(\widetilde{V}_{K},  V_{K}\right):=\left\|\widetilde{V}_{K} \widetilde{V}_{K}^{T}-V_{K} V_{K}^{T}\right\|_{F}$, it is decomposed into the sample variance term $\rho\left(\widetilde{V}_{K},  V_{K}^{*}\right)$ and the bias term $\rho\left(V_{K}^{*},  V_{K}\right)$,
where $V_{K}^{*}=\left(v_{1}^{*},  \ldots,  v_{K}^{*}\right) \in \mathbb{R}^{d \times K} $ is defined as the top $K$ eigenvectors of $\Sigma^{*}:=\frac{1}{m}\sum_{\ell=1}^{m}\mathbb{E}\left(\widehat{V}_{K}^{(\ell)} \widehat{V}_{K}^{(\ell) T}\right)$ and $V_{K}$ is the top $K$ eigenvectors of $\Sigma^{(\ell)}$. For further details on the statistical error, the reader is referred to  Section 4 of Fan et al. (2019)\cite{fan2019distributed}.
%------------------------------------------------

\section{\textbf{Main results}} \label{SecMain}
\subsection{\textbf{Robust covariance estimation}} \label{Robust covariance estimation}
This section initially presents the definitions of the robust covariance matrix estimators of Minsker (2018)\cite{Minsker_2018} and Ke et al. (2019)\cite{User-friendly}.\par

Minsker's shrinkage covariance estimator is defined as follows: Assume that $n$ i.i.d. $d$ dimensional random vectors $\left\{X_{i}\right\}_{i=1}^{n}$ have finite  $v^{\alpha} :=\left\|\mathbb{E} |X_{i}X_{i}^{T}|^{\alpha}\right\|_{2}$ for some $\alpha \in(1,2]$. Let
$$  \widehat{\Sigma}_{n}(\alpha,\theta):=\frac{1}{n \theta} \sum_{i=1}^{n} \psi_{\alpha}\left(\theta X_{i} {X_{i}}^{T}\right)=\frac{1}{n \theta} \sum_{i=1}^{n} \psi_{\alpha}\left(\theta\left\|X_{i}\right\|_{2}^{2}\right) \frac{X_{i} X_{i}^{T}}{\left\|X_{i}\right\|_{2}^{2}}$$
where $\psi_{\alpha}(x): \mathbb{R} \mapsto \mathbb{R}$ be a nondecreasing function such that for all $x \in \mathbb{R}$ and $c_{\alpha}=\frac{\alpha-1}{\alpha} \vee \sqrt{\frac{2-\alpha}{\alpha}}$,
$$
-\log \left(1-x+c_{\alpha}|x|^{\alpha}\right) \leq \psi_{\alpha}(x) \leq \log \left(1+x+c_{\alpha}|x|^{\alpha}\right).
$$
\par

Similar to Minsker's estimator, Ke's truncated covariance estimator is constructed as follows:
$$ \widehat{\Sigma}_{n}(\alpha,\tau):=\frac{1}{n} \sum_{i=1}^{n} \psi_{\tau}\left(\left\|X_{i}\right\|_{2}^{2}\right) \frac{X_{i} X_{i}^{T}}{\left\|X_{i}\right\|_{2}^{2}},$$
where $\psi_{\tau}(x)=(|x| \wedge \tau) \operatorname{sign}(x)$.\par
The tuning parameters $\theta$ and $\tau$ control the bias-robustness trade-off. Because of the above concise expression, we utilize these estimators in Algorithm \ref{Al} such that the computational efficiency is comparable with that of the sample covariance matrix.\par

In Lemma \ref{lemma}, we slightly ameliorate Theorem 3.3 and 3.2 of Minsker (2018)\cite{Minsker_2018} for the Minsker's shrinkage estimator.
\begin{lemma}\label{lemma} Let $X \in \mathbb{R}^{d}$ be a random vector such that $\mathbb{E} X=0_{d \times 1}$,  for some $\alpha \in(1,2]$,  $\left\|\mathbb{E}|X X^{T}|^{\alpha}\right\|_{2}<\infty$ and $\Sigma=\mathbb{E}\left[XX^{T}\right]$. Assume that $X_{1},  \ldots,  X_{n}$ are i. i. d.  copies of $X$ and
$$
v^{\alpha} :=\left\|\mathbb{E}|X_{i} X_{i}^{T}|^{\alpha}\right\|_{2},\;\; \bar{d}_{\alpha}:=\frac{\operatorname{tr}\left(\mathbb{E} |X_{i}X_{i}^{T}|^{\alpha}\right)}{\left\|\mathbb{E} |X_{i}X_{i}^{T}|^{\alpha}\right\|_{2}}.$$
For $\theta>0$, the following inequality holds:
\[
\operatorname{P}\left(\left\|\widehat{\Sigma}_{n}(\alpha,\theta)-\Sigma\right\|_{2} \geq t\right) \leq 2 \bar{d}_{\alpha}e^{-\theta nt}\left(e^{c_{\alpha}n \theta^{\alpha}v^{\alpha}}+c_{\alpha}n \theta^{\alpha}v^{\alpha}\right)\left(1+\frac{2}{\theta nt}+\frac{2}{(\theta nt)^2}\right).
\]
\end{lemma}

The following lemma generalizes Ke's truncated estimator to the case of bounded $(2+\epsilon)$-th moment, with the proof being similar to Lemma \ref{lemma}. The difference is that $\psi_{\tau}(x)=\tau \psi_{1}(x / \tau)$ and the inequality that$-\log \left(1-x+|x|^{\alpha}\right) \leq  \psi_{1}(x)
\leq \log \left(1+x+|x|^{\alpha}\right)$ $\operatorname{for} \operatorname{all} x\in \mathbb{R},\alpha \in (1,2]$ are used. Therefore, the corresponding details are omitted.
\begin{lemma}\label{lemma1.5}
Under the condition of Lemma \ref{lemma}, for any positive $t$ and $\tau$, we have
$$\operatorname{P}\left(\left\|\widehat{\Sigma}_{n}(\alpha,\tau)-\Sigma\right\|_{2} \geq t\right) \leq 2\bar{d}_{\alpha}  e^{-nt/\tau}\left(e^{n\tau^{-\alpha} {\sigma}^{\alpha}}+n\tau^{-\alpha} {\sigma}^{\alpha}\right)\left(1+\frac{2\tau}{nt}+\frac{2\tau^{2}}{(nt)^{2}}\right)$$
where ${\sigma}^{\alpha} :=\left\|\mathbb{E}|X_{i} X_{i}^{T}|^{\alpha}\right\|_{2}$.
\end{lemma}
Next, we further study the deviation bound of Ke's element-wise truncated covariance estimator under max norm, when the samples have only finite $(2+\epsilon)$-th moment. For $1\leq k, s \leq d$, denote $\widehat{\sigma}_{k,s}:=\frac{1}{n}\sum_{i=1}^{n}\psi_{\tau_{k,s}}\left({X_{i}}_{(k)}{X_{i}}_{(s)}\right).$ We weaken the moment assumption of the above robust covariance estimator to show the following theorem.
\begin{theorem}\label{theorem1}
Let $X_{1},  \ldots,  X_{n}$ be $n$ i. i. d. random vectors with $\mathbb{E} X_{1}=0_{d \times 1},  \Sigma=\mathbb{E}X_{1}X_{1}^{T}$. Assume $M:=\max_{k,s \in [d]}\sqrt[\alpha]{{\mathbb{E}\left|{X_{i}}_{(k)}{X_{i}}_{(s)}\right|^{\alpha}}}<\infty$ for some $\alpha \in (1, 2]$, and $\widehat{\Sigma}_{(k,s)}=\widehat{\sigma}_{k,s}$. With probability at least $1-(1+d^{-1})\delta$, we have
$$
\left\|\widehat{\Sigma}-\Sigma\right\|_{\max}\leq 2 M \left(\frac{2\log d + \log \delta^{-1}}{n}\right)^{\frac{\alpha-1}{\alpha}}
$$
where $\tau_{k,s}=\sqrt[\alpha]{n \mathbb{E}\left|{X_{i}}_{(k)}{X_{i}}_{(s)}\right|^{\alpha}/(2\log d-\log \delta)}$.
\end{theorem}

In the detailed numerical simulation study of Ke et al. (2019)\cite{User-friendly}, they demonstrated the statistical performance superiority of the truncated covariance estimator. Although the theoretical properties of their estimator hold only when the fourth moment exists, it still possesses excellent simulation results when the random samples follow $t$ distribution with three degrees of freedom. We believe that our generalized truncated estimator has a better effect than that of them since the selection of our tuning parameter $\tau$ conforms better to the $(2+\epsilon)$-th moment property of random samples.\par

\subsection{\textbf{Distributed estimation with finite moment}} \label{Distributed}
In this subsection, we utilize the above-truncated covariance estimator to analyze the statistical error rate of the resulting distributed estimator $\widetilde{V}_{K}$. We omit to present similar results using the shrinkage covariance estimator to avoid redundancy. We slightly modify Theorem 2 in Fan et al. (2019)\cite{fan2019distributed} to elucidate that when random samples have symmetric innovation, the bias term $\rho\left(V_{K}^{*},  V_{K}\right)$ is equal to $0$ laying the foundation for the proof of Theorem \ref{theorem5}.

\begin{lemma}\label{lemma4} Assume that in Algorithm \ref{Al}, the $\ell$-th machine has $n$ i. i. d. samples $\left\{X_{i}^{(\ell)}\right\}_{i=1}^{n}$ which satisfies the condition of Lemma \ref{lemma}. If $\left\{X_{i}^{(\ell)}\right\}_{i=1}^{n}$ are symmetric innovation and $\left\|\mathbb{E}\widehat{{V}}_{K}^{(\ell)} \widehat{{V}}_{K}^{(\ell) T}-V_{K} V_{K}^{T}\right\|_{2}<1 / 2, $ we have $\rho\left(V_{K}^{*{(\ell)}},  V_{K}\right)=0$ where $V_{K}^{*(\ell)}$ is denoted as the top $K$ eigenspaces of $\mathbb{E}\widehat{{V}}_{K}^{(\ell)} \widehat{{V}}_{K}^{(\ell) T}$.
\end{lemma}

After finishing Lemma \ref{lemma4}, we combine the above results to derive the statistical error rate of $\widetilde{V}_{K}$ when the total $N$ samples have symmetric innovation and the same top $K$ principal eigenspaces $V_{K}$.

\begin{theorem}\label{theorem5} In Algorithm \ref{Al}, suppose that $N$ random samples satisfying the symmetric innovation and the condition of Lemma \ref{lemma} with $\alpha=2$, are distributed on $m$ machines. In each machine, $n$ samples are i.i.d., but between different machines, samples are only independent and have different tails as well as various covariance structures $\{\Sigma^{(\ell)}\}_{i=1}^{m}$. If we choose $\tau_{(\ell)}=O\left(\sigma_{(\ell)} \cdot \sqrt{n}\right)$ for $\ell \in [m]$, there exist constants $C_{1}$ and $C_{2}$ such that as long as $n \geq C_{1} \max_{\ell \in [m]}\left(\bar{d}_{(\ell)}\frac{\sigma_{(\ell)}}{\Delta_{(\ell)}}\right)^2K$,
\[
\left\|\rho\left(\widetilde{{V}}_{K},  {V}_{K}\right)\right\|_{\psi_{1}} \leq C_{2}\sqrt{\frac{1}{m}\sum_{\ell=1}^{m}\left(\bar{d}_{(\ell)} \frac{\sigma_{(\ell)}}{\Delta_{(\ell)}}\right)^{2}} \sqrt{\frac{K}{N}},
\]
where $\Delta_{(\ell)}=\lambda_{K}(\Sigma^{(\ell)})-\lambda_{K+1}(\Sigma^{(\ell)})>0$.
\end{theorem}
\begin{remark}\label{Remark3}
$\bar{d}_{(\ell)}$ in the upper bound is indispensable and often much smaller than $d$. In the simulation study presented next, we will illustrate that under some mild conditions, the upper bound of the truncated covariance estimator with $\psi_{1}$-norm is consistent with the statistical error rate $O_{P}(\sqrt{d \log d/n})$ of $\ell_{4}$-norm truncated sample covariance matrix constructed by Fan, Wang and Zhu (2021)\cite{A shrinkage principle}.
\end{remark}

\begin{remark}\label{Remark4}
When the sample follows the elliptical distribution without finite kurtosis, the multivariate Kendall's tau statistic (Han and Liu (2018)\cite{ECA}) can be employed in Algorithm \ref{Al} such that the statistical error rate of estimating eigenspaces is similar to that of the sub-Gaussian tail. Furthermore, if the magnitude of $K$ is in practice unknown, the weighted distributed PCA method proposed by Bhaskara and Wijewardena (2019)\cite{distributed averaging} can reinforce the practicability.
\end{remark}

Next, we consider the asymmetric case with finite $6$-th moment.

\begin{lemma}\label{lemma4.5} Assume that $n$ i.i.d. samples $\left\{X_{i}\right\}_{i=1}^{n}$ on the $1$-th machine satisfy $\mathbb{E}\left(v^{T} X_{i}\right)^{6} \leq R_{(1)}^{\prime}$  for $\forall v \in \mathcal{S}^{d-1}$, $\mathbb{E}X_{i}=0_{d\times1}$ and $\Sigma^{(1)}=\mathbb{E}X_{i}X_{i}^{T}$. Then there are constants $C_{1}$ and $C_{2}$ such that $\left\|\mathbb{E}\widehat{{V}}_{K}^{(1)} \widehat{{V}}_{K}^{(1) T}-{V}_{K} {V}_{K}^{T}\right\|_{F} \leq C_{1} \left(\bar{d}_{(1)}\frac{\sigma_{(1)}}{\Delta_{(1)}}\right)^2\frac{\sqrt{K}}{n}+C_{2}
\frac{R_{(1)}^{\prime}d^{2}}{\sigma_{(1)}^2\Delta_{(1)}}\frac{\sqrt{K}}{n}.$
\end{lemma}
Based on Lemma \ref{lemma4.5}, we have the following theorem.
\begin{theorem}\label{theorem6}Suppose that $N$ random samples that are distributed on $m$ machines, satisfy the condition of Theorem \ref{theorem5} and Lemma \ref{lemma4.5} except for the symmetric innovation. There exist some constants $\{C_{i}\}_{i=1}^3$ such that
\[
\left\|\rho\left(\widetilde{{V}}_{K},  {V}_{K}\right)\right\|_{\psi_{1}}\leq C_{3}\sqrt{\frac{1}{m}\sum_{\ell=1}^{m}\left(\bar{d}_{(\ell)} \frac{\sigma_{(\ell)}}{\Delta_{(\ell)}}\right)^{2}} \sqrt{\frac{K}{N}},
\]
as long as $n\geq C_{1}\max\left(\frac{1}{m}\sum_{\ell=1}^{m}\left(\bar{d}_{(\ell)} \frac{\sigma_{(\ell)}}{\Delta_{(\ell)}}\right)^{2},\frac{1}{m}\left(\sum_{\ell=1}^{m}\frac{R_{(\ell)}^{\prime}d^{2}}
{\sigma_{(\ell)}^2\Delta_{(\ell)}}\right)^2/\sum_{\ell=1}^{m}\left(\bar{d}_{(\ell)} \frac{\sigma_{(\ell)}}{\Delta_{(\ell)}}\right)^{2}\right)m$ and $n \geq C_{2} \max_{\ell \in [m]}\left(\bar{d}_{(\ell)}\frac{\sigma_{(\ell)}}{\Delta_{(\ell)}}\right)^2K$.\par
\end{theorem}
\begin{remark}\label{Remark5}
When the total $N$ samples are i.i.d., Theorem \ref{theorem6} shows that by adding the additional requirement of $n\gtrsim\max\left(\left(\bar{d}\frac{\sigma }{\Delta}\right)^2,\left(\frac{R^{\prime}d^{2}}{\sigma^3 \bar{d}}\right)^2\right)m$, the robust distributed PCA method possesses the same statistical error rate as the one established in Theorem \ref{theorem5}.
\end{remark}

\section{\textbf{Numerical experiment}}\label{SecSimulation}
In this section, we first follow the simulation framework of Fan et al. (2019)\cite{fan2019distributed} to validate that the statistical error rate of $\widetilde{{V}}_{K}$ differs from the one of Gaussian tails when the samples are symmetric innovation and have only polynomial and exponential-type tails. Then we compare the statistical performance of $\widetilde{{V}}_{K}$ produced by Algorithm \ref{Al} with the distributed PCA algorithm proposed by Fan et al. (2019)\cite{fan2019distributed} under the finite moment assumption and outlier-robust distributed setting. The major conclusion is that the proposed robust distributed estimator performs better than Fan's estimator.
\subsection{\textbf{Leading eigenvectors of spiked covariance matrix}}\label{subsect1}
To validate our results, we consider two types of distributions with fourth bounded moments:\par
1. Polynomial-type tail case: We assume that $\left\{X_{i}\right\}_{i=1}^{N}$ i. i. d.  follow $T_{d}(5, 0_{d\times1}, \Sigma), $ where $\Sigma=\operatorname{diag}(\lambda,  \lambda / 2,  \lambda / 4, 1,  \ldots,  1)$ . The population covariance matrix is $5/3\operatorname{diag}(\lambda, \lambda / 2,\lambda / 4, 1,  \ldots, 1).$\par
2. Exponential-type tail case: Let $X_{(1)}, X_{(2)}$ and $X_{(3)}$ follow Laplace distributions $\text{La}(0,\sqrt{\lambda})$, $\text{La}(0,\sqrt{\lambda/2})$ and $\text{La}(0,\sqrt{\lambda/4})$, respectively, and for $j\in \{4, 5, \ldots, d\}, X_{(j)}$ follows $\text{La}(0, 1)$. Denote $X=\left(X_{(1)},  \ldots,  X_{(d)}\right)^{T}$. We consider that $\left\{X_{i}\right\}_{i=1}^{N}$ i. i. d.  follow the above distribution. The population covariance matrix is $2\operatorname{diag}(\lambda,  \lambda / 2,  \lambda / 4, 1,  \ldots,  1)$.\par
Lemma 2.3 in  Minsker and Wei (2017)\cite{10.5555/3294996.3295045} shows that
$r\left(\Sigma\right)\left\|\Sigma\right\|_{2}^{2} \leq v^{2} \leq R^{2} r\left(\Sigma\right)\left\|\Sigma\right\|_{2}^{2}.$ Hence, by Theorem \ref{theorem5}, we can easily obtain that when $\lambda=O(1)$,
\begin{equation}\label{eq13}
\left\|\left\|\widehat{\Sigma}_{n}(2,\theta)-\Sigma\right\|_{2}\right\|_{\psi_{1}} \lesssim \bar{d}{\frac{v}{\sqrt{n}}}\lesssim \bar{d}R{\frac{ \sqrt{r\left(\Sigma\right)}\left\|\Sigma\right\|_{2}}{\sqrt{n}}}=\bar{d}R\lambda
\sqrt{\frac{r}{n}}=O\left(\bar{d}\sqrt{\frac{d}{n}}\right).
\end{equation}
When $\lambda=O(d)$,
\begin{equation}\label{eq14}
\left\|\rho\left(\widetilde{V}_{K}, V_{K}\right)\right\|_{\psi_{1}}=O\left(\frac{\bar{d}R\|\Sigma\|_{2}}{\lambda_{K}-\lambda_{K+1}} \sqrt{\frac{K r(\Sigma)}{N}}\right)=O\left(\bar{d}\sqrt{\frac{d}{m n \delta}}\right)
\end{equation}
where $\delta:=\lambda_{K}-\lambda_{K+1}.$ On the one hand, as Fan, Wang and Zhu (2021)\cite{A shrinkage principle} mentioned, $R$ can vary with the increase of the dimension $d$ or other parameters. We hide $R$ in $\bar{d}$ and do not regard it as a single independent variable. On the other hand, the following simulation results show that the effective dimension $\bar{d}$ in (\ref{eq14}) has a minor effect on the statistical error rate. To confirm the result of Theorem \ref{theorem5}, we keep $\theta=O\left(\frac{1}{v \sqrt{n}}\right)$ where $v$ is estimated via the empirical moment estimator $\sqrt{\|\frac{1}{n}\sum_{i=1}^{n}\|X_{i}\|_{2}X_{i}X_{i}^{T}\|_{2}}$ and the constant is controlled by cross-validation.\par

\begin{figure}[H]
  \centering
  \includegraphics[width=6in]{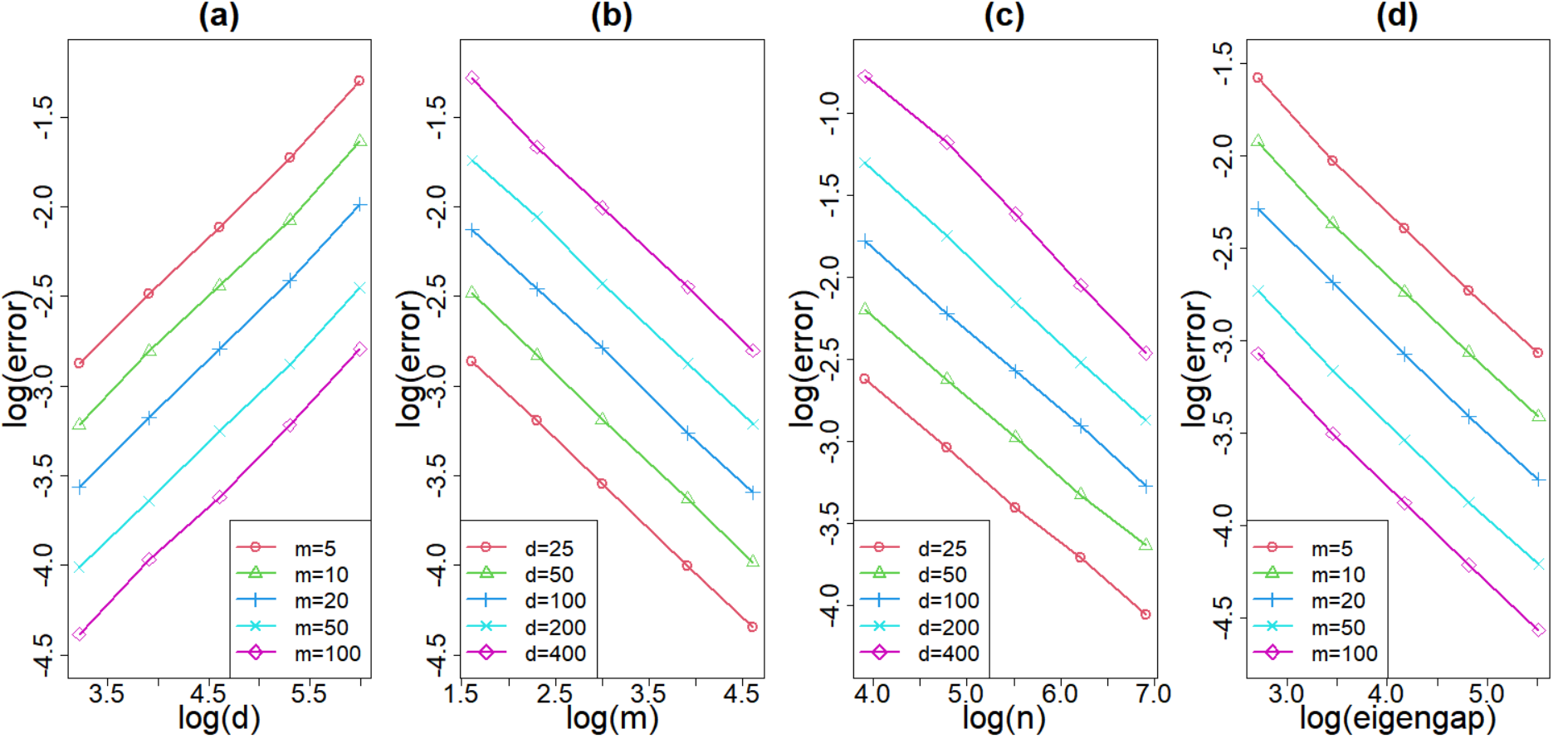}
  \caption{Polynomial-type tail case: The statistical error rate with respect to:
$(a)$ the dimension $d$ when $\lambda=50$ and $n=1000$; (b) the number of servers $m$ when $\lambda=50$ and $n=1000$; (c) the subsample size $n$ when $\lambda=50$ and $m=50$; (d) the eigengap $\delta$ when $d=200$ and $n=1000$. }
\label{fig1}
\end{figure}
\begin{figure}[H]
  \centering
  \includegraphics[width=6in]{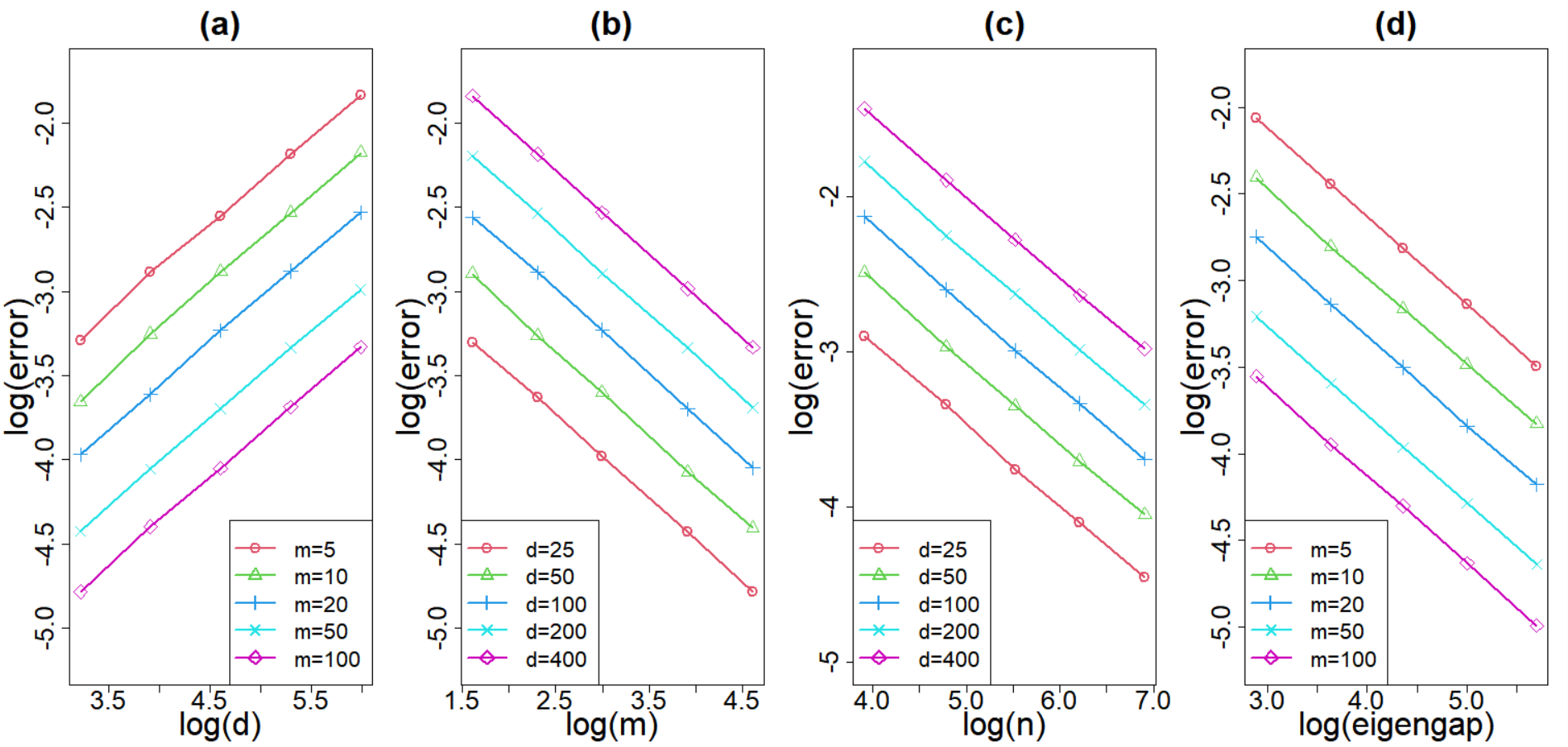}
  \caption{Exponential-type tail case: The statistical error rate with respect to:
(a) the dimension $d$ when $\lambda=50$ and $n=1000$; (b) the number of servers $m$ when $\lambda=50$ and $n=1000$; (c) the subsample size $n$ when $\lambda=50$ and $m=50$; (d) the eigengap $\delta$ when $d=200$ and $n=1000$. }
\label{fig2}
\end{figure}
Both plots above depict the marginal relationship between $\rho\left(\widetilde{V}_{K},  V_{K}\right)$ and the four parameters $d$, $m$,  $n$ and $\delta.$ Each data point on the plots is based on 50 independent Monte Carlo simulations. To highlight the difference against the result of Fan et al. (2019)\cite{fan2019distributed}, we fit the regression model based on all the data points $\left(d,  m,  n,  \delta,  \rho\left(\widetilde{V}_{K},  V_{K}\right)\right)$ from the four plots in Figure \ref{fig1} and \ref{fig2}, respectively.  The fitting results are that for polynomial-type tail case,  $\widehat{\beta}_{d}=0. 5790,  \widehat{\beta}_{m}=-0. 4987,  \widehat{\beta}_{n}=-0. 5085$ and $\widehat{\beta}_{\delta}=-0. 5420$ with the multiple $R^{2}=0. 9975, $ and for the exponential-type tail case, $\widehat{\beta}_{d}=0. 5242,  \widehat{\beta}_{m}=-0. 4997,  \widehat{\beta}_{n}=-0. 5172$ and $\widehat{\beta}_{\delta}=-0. 5079$ with the multiple $R^{2}=0. 9998.$ Note that in both cases $\widehat{\beta}_{d}$ exceed 0.5 which is consistent with (\ref{eq14}). This authenticates that $\bar{d}=O(\sqrt{\log d})$ is much smaller than $\sqrt{d}$ and implies that the convergence rate in (\ref{eq13}) is the same as the error rate $O_{P} (\sqrt{d\log d/n})$ of Fan, Wang and Zhu (2021)\cite{A shrinkage principle}. $\widehat{\beta}_{m}$ and $\widehat{\beta}_{n}$ are all consistent with (\ref{eq14}) and indicate that $\rho\left(\widetilde{V}_{K},  V_{K}\right)$ is proportional to $m^{-\frac{1}{2}}$ and $n^{-\frac{1}{2}}$, respectively. Therefore, the robust distributed estimator achieves nearly the same statistical convergence rate as that of sub-Gaussian samples, up to just a small parameter $\bar{d}$.\par

\subsection{\textbf{Comparison under finite moment assumption}}\label{subsect4}
We denote the robust distributed PCA of this paper by RDP and the distributed PCA proposed by Fan et al. (2019)\cite{fan2019distributed} by DP, respectively. To investigate the robustness against the low order moments of samples and demonstrate the superiority of our method, we offer some numerical results in three scenarios described below to compare the statistical performance between RDP and DP when the truncated covariance estimator is applied.\par

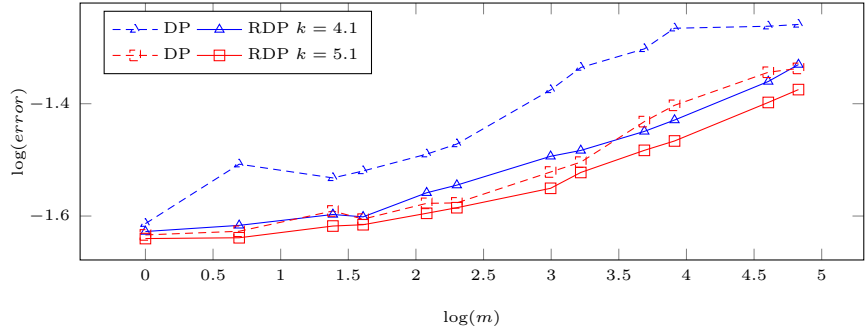
\begin{figure}[H]
\centering
\subfigure{
\begin{tikzpicture}
\pgfplotsset{every axis legend/.append style={at={(0.5,0.5)},anchor=south},every axis y label/.append style={at={(0.04,0.5)}}}
\begin{axis}[xlabel = $\log(m)$,ylabel = {$\log(error)$}, legend columns=2,legend pos=north west,title style={font=\tiny},legend style={font=\tiny},font=\tiny,height=5cm, width=12cm]
\addplot[
    densely dashed,
    color=blue,
    mark=triangle
    ]
    coordinates {
(0,-1.613557)(0.6931,-1.50767)(1.3862,-1.532091)(1.6094,-1.51977)(2.0794,-1.489812)(2.3025,-1.471757)(2.9957,-1.375469)(3.2188,-1.335563)
(3.6888,-1.302345)(3.9120,-1.265587)(4.6051,-1.262023)(4.8283,-1.258834)};
\addplot[
    color=blue,
    mark=triangle]
    coordinates {
(0,-1.627551)(0.6931,-1.616610)(1.3862,-1.597207)(1.6094,-1.601454)(2.0794,-1.558569)(2.3025,-1.544947)(2.9957,-1.493617)(3.2188,-1.483202)
(3.6888,-1.44928)(3.9120,-1.429062)(4.6051,-1.360404)(4.8283,-1.330269) };
\addplot[
    densely dashed,
    color=red,
    mark=square
    ]
    coordinates {
(0,-1.633624)(0.6931,-1.627111)(1.3862,-1.590968)(1.6094,-1.605529)(2.0794,-1.577407)(2.3025,-1.576591)(2.9957,-1.521572)(3.2188,-1.504016)
(3.6888,-1.431840)(3.9120,-1.402845)(4.6051,-1.343968)(4.8283,-1.337094) };
\addplot[
    color=red,
    mark=square]
    coordinates {
(0,-1.640109)(0.6931,-1.638682)(1.3862,-1.617789)(1.6094,-1.61548)(2.0794,-1.595095)(2.3025,-1.585161)(2.9957,-1.550602)(3.2188,-1.522404)
(3.6888,-1.483018)(3.9120,-1.466323)(4.6051,-1.397738)(4.8283,-1.374843) };
\legend{DP, RDP $k=4.1$, DP, RDP $k=5.1$}
\end{axis}
\end{tikzpicture}}
\caption{The trend of the estimation error with varying $m$, when $\lambda=50$ and $d=600$.}\label{fig3}
\end{figure}

In the first scenario, we consider that the samples follow Pareto distributions with $\Sigma=\operatorname{diag}(\lambda, \lambda / 2,\lambda / 4, 1,  \ldots, 1)$ and different shape parameters $k \in\{4.1, 5.1\}$. The total sample size is $N=5000$. Figure \ref{fig3} illustrates that after carefully choosing $\tau$ and $50$ independent Monte Carlo simulations, the average estimation error slightly rises as $m$ increases.

For each of the following two scenarios, we employ multivariate $t$ distributions with different degrees of freedom $\nu \in \{4.1, 4.5, 5.0, 5.5, 6.0\}$ and calculate $\rho\left(\widetilde{V}_{K},  V_{K}\right)$  based on 50 independent Monte Carlo simulations. Note that if we set the tuning parameter $\tau=O\left(\sqrt{\|\mathbb{E}|XX^{T}|^{2}\|_{2}} \sqrt{n}\right)$, the constant to be adjusted, which determines the performance of RDP, is hard to control. Therefore, we solve the following equation of Ke et al. (2019)\cite{User-friendly} to adaptively obtain $\tau$:
$$\left\|\frac{1}{\tau^{2}} \sum_{i=1}^{n}\left(\left\|X_{i}\right\|_{2}^{2} \bigwedge \tau\right)^{2} \frac{X_{i} X_{i}^{T}}{\left\|X_{i}\right\|_{2}^{2}}\right\|_{2}=\log(2d)+\log(n).$$
In Table \ref{tab1}, we generate $n \in \{200, 400, 800, 1600, 2000\}$ i.i.d. subsamples to obtain the logarithmic error when $m=20, \lambda=50$ and $d=100, 200, 300, 400$. In Table \ref{tab2}, we increase the number of servers $m$ from 5 to 100 to compare the two methods when $n=400, d=200$ and $\lambda=10, 20, 40, 80$. The experimental results indicate that when $\nu$ and $\lambda$ decrease, the estimation result turns to be not viable for the DP method.\par

From Figure \ref{fig3} and Table 1-2, we conclude that RDP provides the following evidence-based advantages:
\par
1.  When all parameters (subsample size $n$, dimension $d$, number of servers $m$, $\lambda$ and degree of freedom $\nu$) are fixed, the statistical performance of RDP is superior to DP;\par
2.  For the asymmetric distribution with only finite $4$-th moment, the result of Theorem \ref{theorem6} still holds;\par
3.  The proposed algorithm is more appealing when $d/n$ increases or $\lambda$ is small.

\begin{table}[H]
    \centering
    \fontsize{9}{10}\selectfont
    \caption{Comparison between RDP and DP with respect to the subsample size $n$ under multivariate t distribution with different degrees of freedom $\nu$ when $m=20, \lambda=50$. $\widehat{\beta}_{n}$ represents the regression coefficient between the logarithmic error and $\log(n)$.}
    \label{tab1}
    \begin{tabular}{ccccccccccc}
    \toprule
    \multirow{2}{*}{$\nu$}&\multirow{2}{*}{$n$}&\multicolumn{2}{c}{$d=100$}& \multicolumn{2}{c}{$d=200$}&\multicolumn{2}{c}{$d=300$} &\multicolumn{2}{c}{$d=400$}\cr
   \cmidrule(lr){3-4}\cmidrule {5-6}\cmidrule(lr){7-8}\cmidrule {9-10}

  &  & RDP & DP&RDP &DP&RDP& DP & RDP&DP\cr
    \midrule
      $6.0$ &200 & -2.2325 & -2.0827&-1.9269 &-1.6554& -1.7566 & -1.3713& -1.6424 & -1.1169
  \\
       &400 & -2.5650 & -2.4228& -2.2294&-2.0292& -2.0581 &-1.7650 & -1.9329 &-1.5615
  \\
       &800 &-2.8737 & -2.7645& -2.5413 & -2.3969& -2.3528 & -2.1512 &-2.2227 & -1.9573
  \\
       &1600 & -3.2059 &-3.1420 & -2.8584 &-2.7586& -2.6712 &-2.5187 & -2.5309 &-2.3470
  \\
       &2000 & -3.3016& -3.2430 & -2.9683 & -2.8764 & -2.7689 & -2.6333 & -2.6345 &-2.4646
  \\ \midrule
      \multicolumn{2}{c}{$\widehat{\beta}_{n}$}  & -0.4641 & -0.5075 &-0.4518&-0.5296& -0.4400 & -0.5480&  -0.4301 & -0.5825
  \\  \midrule
      $5.5$ &200 & -2.2286  & -2.0251&-1.9317 &-1.5813& -1.7523 & -1.2800& -1.6416 & -1.0180
  \\
       &400 & -2.5328 & -2.3769 & -2.2156 &-1.9703&-2.0381 &-1.7031 & -1.9210 &-1.4461
  \\
       &800 &-2.8626 & -2.7279 & -2.5229 & -2.3176&-2.3374 & -2.0669 & -2.2080 &-1.8705
  \\
       &1600 & -3.1678 &-3.0486 & -2.8243 &-2.6713& -2.6473 &-2.4526& -2.5104 &-2.2730
  \\
       &2000 & -3.2651  &-3.1729  & -2.9329  &-2.8019 & -2.7427  &-2.5604 & -2.6144& -2.3738
  \\ \midrule
      \multicolumn{2}{c}{$\widehat{\beta}_{n}$}  & -0.4523 & -0.4950 & -0.4349 & -0.5240 & -0.4320 & -0.5538 & -0.4221 &-0.5930
  \\ \midrule
   $5.0$&200 & -2.2110 & -1.9590 & -1.9181 &-1.4856& -1.7480  &-1.1718 & -1.6394 &-0.9013
  \\
       &400 & -2.5166& -2.2953 & -2.2055 &-1.8712& -2.0304 &-1.5734 & -1.9089 &-1.3314
  \\
       &800 &-2.8225 & -2.6353 & -2.4948 &-2.2124 & -2.3141 &-1.9370 & -2.1944 &-1.7359
  \\
       &1600 & -3.1350& -2.9486 & -2.7995 &-2.5627& -2.6105 &-2.3193 & -2.4877 &-2.1218
  \\
       &2000 & -3.2322 & -3.0773 &-2.8940 & -2.6714 & -2.7123& -2.4527 & -2.5740 & -2.2704
  \\ \midrule
      \multicolumn{2}{c}{$\widehat{\beta}_{n}$} & -0.4441 & -0.4815 & -0.4246&-0.5120 & -0.4179 &-0.5515 & -0.4085 & -0.5885
  \\ \midrule
  $4.5$ &200 & -2.1925  & -1.8714 & -1.8972 &-1.3873& -1.7507 &-1.0628 & -1.6360 &-0.7544
  \\
       &400 &-2.4865 & -2.1903 & -2.1869 &-1.7408& -2.0133 &-1.4352 & -1.8996 &-1.1622
  \\
       &800 & -2.7907 &-2.5359 &-2.4576 &-2.0955 & -2.2896 &-1.8213 &-2.1728 &-1.5605
  \\
       &1600 & -3.0743 &-2.8379 & -2.7627 &-2.4133& -2.5814 &-2.1581& -2.4616 &-1.9716
  \\
       &2000 & -3.1807 & -2.9394 & -2.8643 & -2.5613& -2.6727 & -2.2653& -2.5479 & -2.0959
  \\ \midrule
      \multicolumn{2}{c}{$\widehat{\beta}_{n}$}  & -0.4276 & -0.4652 & -0.4180& -0.5026& -0.4020 & -0.5234 & -0.3979 & -0.5831
  \\ \midrule
      $4.1$ &200& -2.1796 &-1.7943 & -1.8958 &-1.2879& -1.7413 &-0.9472 & -1.6291 &-0.6703
  \\
       &400 & -2.4849 & -2.0848 & -2.1678 &-1.6312& -2.0016 &-1.2930 & -1.8873 &-0.9888
  \\
       &800 & -2.7613 & -2.3854 & -2.4381 & -1.9842 & -2.2654 & -1.6539 & -2.1512 &-1.4375
  \\
       &1600 & -3.0578 & -2.6833 & -2.7275 &-2.2743& -2.5530 &-1.9838 & -2.4303 &-1.7689
  \\
       &2000 & -3.1427 & -2.7919 & -2.8216 &-2.3681 & -2.6459 &-2.0903 & -2.5210 &-1.9035
  \\ \midrule
      \multicolumn{2}{c}{$\widehat{\beta}_{n}$}  & -0.4176 & -0.4320 &-0.4018 &-0.4697 & -0.3932  & -0.4976 & -0.3877 & -0.5416
  \\
    \bottomrule
\end{tabular}\vspace{0cm}
\end{table}
\begin{table}[H]
    \centering
    \fontsize{9}{10}\selectfont
    \caption{Comparison between RDP and DP with respect to the number of servers $m$ under multivariate t distribution with different degrees of freedom $\nu$ when $n=400, d=200$. Numbers in bold denote the estimation error greater than 0. $\widehat{\beta}_{m}$ represents the regression coefficient between the logarithmic error and $\log(m)$.}
    \label{tab2}
    \begin{tabular}{ccccccccccc}
    \toprule
    \multirow{2}{*}{$\nu$}&\multirow{2}{*}{$m$}&
   \multicolumn{2}{c}{$\lambda=10$}&\multicolumn{2}{c}{$\lambda=20$}& \multicolumn{2}{c}{$\lambda=40$}&\multicolumn{2}{c}{$\lambda=80$}\cr
   \cmidrule(lr){3-4}\cmidrule {5-6}\cmidrule(lr){7-8}\cmidrule {9-10}

  &  & RDP & DP&RDP &DP&RDP& DP & RDP&DP\cr
  \midrule
       $6.0$ &5 & -0.4768&\textbf{0.3260}&-1.0391&-0.4562&-1.4266 &-1.1926 &-1.7849 &-1.6436
  \\
       &10  & -0.7778&\textbf{0.2021}&-1.3758&-0.8740&-1.7711 & -1.5265 & -2.1316&-1.9811
  \\
       &25 & -1.2141 & -0.0521&-1.8191&-1.2951&-2.2306 & -1.9934 & -2.5837 &-2.4445
  \\
       &50 & -1.5500 & -0.3331&-2.1558&-1.6384&-2.5659 & -2.3275 & -2.9229 &-2.7748
  \\
       &100 & -1.8943 & -0.6161 &-2.5102&-2.0051&-2.9144 &-2.6629 & -3.2766& -3.1365
  \\ \midrule
      \multicolumn{2}{c}{$\widehat{\beta}_{m}$}  & -0.4744 & -0.3172 &-0.4895 &-0.5075  &-0.4961  &-0.4925 & -0.4965  & -0.4973
  \\  \midrule
      $5.5$ &5 & -0.4562 &\textbf{0.3659}& -1.0185&-0.3527& -1.4175&-1.0785 & -1.7607&-1.5707
  \\
       &10  & -0.7694 &\textbf{0.2742}& -1.3626&-0.7370& -1.7493& -1.4270 & -2.1076&-1.9240
  \\
       &25 & -1.2065 &\textbf{0.0197}& -1.8185&-1.1778& -2.2079& -1.9019 & -2.5631& -2.3892
  \\
       &50 & -1.5533 & -0.2136& -2.1599&-1.5313& -2.5576&-2.2449 & -2.9123 & -2.7144
  \\
       &100 & -1.8906 & -0.5363 & -2.4955&-1.8579& -2.9037&-2.6026 & -3.2554 & -3.0752
  \\ \midrule
      \multicolumn{2}{c}{$\widehat{\beta}_{m}$}  & -0.4804 & -0.3004 &-0.4936 &-0.5005 &-0.4974 &-0.5087 & -0.4992 & -0.4999
  \\ \midrule
      $5.0$ &5 & -0.4543& \textbf{0.3879} &-1.0163 &-0.2432& -1.4058& -0.9906 & -1.7483 & -1.4754
  \\
       &10 &-0.7638& \textbf{0.3044} & -1.3472&-0.5742& -1.7358& -1.3263 & -2.0978 & -1.8462
  \\
       &25 &-1.2023 & \textbf{0.0637} & -1.8103 &-1.0291& -2.1950 &-1.8032 & -2.5531 & -2.2940
  \\
       &50 &-1.5426 & -0.1657 & -2.1507&-1.3518& -2.5406 & -2.1593 & -2.8909 & -2.6528
  \\
       &100 &-1.8781 & -0.4573 & -2.5033&-1.7247& -2.8958 & -2.5079 & -3.2364 & -2.9799
  \\ \midrule
      \multicolumn{2}{c}{$\widehat{\beta}_{m}$}  & -0.4771 & -0.2832 & -0.4970&-0.4920&-0.4979 &-0.5090 & -0.4959 & -0.5020
  \\ \midrule
     $4.5$ &5 & -0.4519& \textbf{0.4228} &-0.9994 &-0.1159& -1.3860 &-0.8647 &-1.7195 &-1.3463
  \\
       &10 & -0.7696 &\textbf{0.3511}&-1.3541 &-0.4304& -1.7114 & -1.2025 & -2.0739 & -1.7388
  \\
       &25 & -1.2020 & \textbf{0.1474} & -1.7870 &-0.8586 & -2.1827 & -1.6835 & -2.5203& -2.1953
  \\
       &50 &-1.5284 & -0.0672 &-2.1456 &-1.1910 & -2.5233 & -2.0196 & -2.8614 & -2.5429
  \\
       &100 & -1.8784 & -0.3802&-2.4907 &-1.5332 & -2.8705 & -2.3740 & -3.2227 & -2.8817
  \\ \midrule
      \multicolumn{2}{c}{$\widehat{\beta}_{m}$}  & -0.4750& -0.2648  &-0.4962 &-0.4728&-0.4976 &-0.5049 &  -0.4989 & -0.5098
  \\ \midrule
      $4.1$ &5 &-0.4465 & \textbf{0.4552} &-0.9919&\textbf{0.0372} &-1.3702 & -0.7142 &-1.6993 &-1.2478
  \\
       &10 & -0.7554& \textbf{0.3676}&-1.3374&-0.2477&-1.7116 & -1.1122 &-2.0424 & -1.6188
  \\
       &25 &-1.1961& \textbf{0.2090} &-1.7852&-0.6790 &-2.1507 & -1.5399& -2.4959 & -2.0849
  \\
       &50 & -1.5309& -0.0266 &-2.1278 &-1.0367 &-2.5002 & -1.8884& -2.8434 &-2.4233
  \\
       &100 & -1.8716& -0.3068 &-2.4780&-1.3691&-2.8556 &-2.2417 & -3.1944&-2.7798
  \\ \midrule
      \multicolumn{2}{c}{$\widehat{\beta}_{m}$}  & -0.4770 & -0.2505 &-0.4949&-0.4737&-0.4943 &-0.5036 & -0.4987  &-0.5090
  \\
    \bottomrule
\end{tabular}\vspace{0cm}
\end{table}

\subsection{\textbf{Outlier-robust experiment and applications}}\label{Byzantine}
\begin{figure}[H]
\centering
\begin{tikzpicture}
\pgfplotsset{every axis legend/.append style={at={(0.5,1)},
anchor=south}, every axis y label/.append style={at={(0.05,0.5)}}}
\begin{axis}[legend columns=1,title style={font=\tiny},xlabel = $c$,
    ylabel = {$\log(error)$},legend style={at={(0.78,0.45)}, anchor=north,legend columns=-1},xtick={2,3,4,5,6}, font=\tiny,height=6cm, width=10cm, xmajorgrids=true]
\addplot[
    color=red,
    mark=triangle]
    coordinates {
    (2,-0.9187757)(2.5,-0.1082227)(3,0.2730619)(3.5,0.3450252)(4,0.4468579)(4.5,0.4933866)(5,0.5681548)};
\addplot[
    color=blue,
    mark=+
    ]
    coordinates {
    (2,-1.165765)(2.5,-1.0014321)(3,-0.7392675)(3.5,-0.4856798)(4,-0.185237981)(4.5,0.01844151)(5,0.003181807)
 };
\addplot[
    color=green,
    mark=square
    ]
    coordinates{
    (2,-1.277008)(2.5,-1.260549)(3,-1.285581)(3.5,-1.254141)(4,-1.253240)(4.5,-1.258041)(5,-1.248971)
};
\legend{DP ,RDP with shrinkage,RDP with truncation}
\end{axis}
\end{tikzpicture}
\caption{Comparison between RDP and DP with varying $c$.}
\label{fig4}
\end{figure}
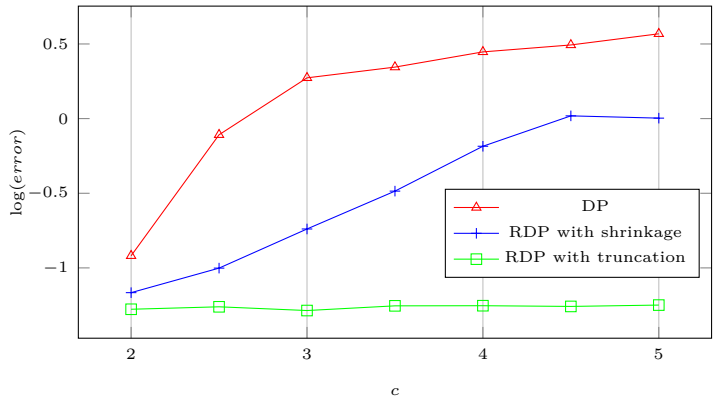

In this subsection, we conduct the experiment to compare the performance of distributed algorithms in the presence of malicious outliers. We partition $N=10^5$ samples into $m=20$ subsamples and the samples i.i.d. follow $d=1000$ dimensional Gaussian distribution with mean zero and $\Sigma=\operatorname{diag}(5, 4, 3, 2, 1, \ldots, 1)$. The outliers are generated as follows: We randomly replace $100$ samples $\left\{X_{i}^{(\ell)}\right\}_{i=1}^{100}$ with $\left\{cX_{i}^{(\ell)}\right\}_{i=1}^{100}$ in each machine where $c\geq 2$ such that the samples are enlarged to moderate values as outliers. The numerical results are presented in Figure \ref{fig4}.
It is evident that as the magnitude of $c$ increases, the proposed algorithm manages a better statistical performance on estimating the top four eigenspaces than the original algorithm, especially on the use of the truncated covariance estimator.\par
Inspired by Chen et al. (2021)\cite{Chen}, we consider the following single index model with finite $(2+\epsilon)$-th moment noise as an application of Algorithm \ref{Al}:
$$
y=f\left(\left\langle\beta^{*}, X\right\rangle\right)+\varepsilon
$$
where $y\in \mathbb{R}$ and $X\in \mathbb{R}^{d}$ are the response and covariate, respectively. $\varepsilon$ which is independent of $X$, is the stochastic error with mean zero and $\mathbb{E}|\varepsilon|^{2+\epsilon}<\infty$. The covariate $X$ follows the $d$-dimensional standard normal distribution. To estimate the direction of the parameter vector $\beta^{\star}$ via second-order Stein's identity, we further assume that the link function $f$ is twice differentiable such that $\mathbb{E}\left[f^{\prime\prime}\left(\left\langle\beta^{*}, X\right\rangle\right)\right]\neq 0$.\par
Under the above assumption, utilizing the second-order Stein's identity (Janzamin et al. (2014)\cite{Janzamin}, Yang et al. (2017)\cite{Yang}), we have
$$
\mathbb{E}\left[y\left(XX^{T}-I_{d}\right)\right]=\mathbb{E}\left[f\left(\left\langle X,\beta^{*}\right\rangle\right)\left(XX^{T}-I_{d}\right)\right]=C \cdot \beta^{*} \beta^{* T}
$$
where $C=2 \mathbb{E}\left[f^{\prime \prime}\left(\left\langle X, \beta^{*}\right\rangle\right)\right]$. Since
$$\left\|\mathbb{E}\left|y\left(XX^{T}-I_{d}\right)\right|^{2}\right\|_{2}\leq
\mathbb{E}\left\|\left|y\left(XX^{T}-I_{d}\right)\right|^{2}\right\|_{2}
=\max\left(\mathbb{E}|y\|X\|_{2}^{2}-y|^{2},\mathbb{E}|y|^{2}\right),$$
by $C_{r}$ inequality and the appropriate choice for $f$, we have $\|\mathbb{E}|yXX^{T}-yI_{d}|^{2}\|_{2}< \infty$.
Therefore, for $n$ i.i.d. sample $\left\{y_{i},X_{i}\right\}_{i=1}^{n}$, we use the random matrix version of Lemma \ref{lemma1.5} to obtain the robust estimator $\widehat{\Sigma}_{n}(\tau):=\frac{1}{n} \sum_{i=1}^{n} \psi_{\tau}\left(y_{i}\left(X_{i}X_{i}^{T}-I_{d}\right)\right)$.\par
We choose $f(x)=x^2$ and $f(x)=x^4/2$ as the link functions in the distributed environment, respectively. Considering that $\varepsilon$ follows $t$ distributions with different degrees of freedom $\nu$ and $\beta^{*}=\beta /\|\beta\|_{2}$ with $\beta \sim N\left(0_{d\times 1}, I_{d}\right)$, then each machine $\ell$ computes the first eigenvector of $\widehat{\Sigma}_{n}^{(\ell)}(\tau)$ and the central machine yields the final estimator $\hat{\beta}$. We use $\rho(\hat{\beta}, \beta^{*}):=\min \left\{\| \|\hat{\beta}\|_{2}^{-1}\hat{\beta}-\beta^{*}\|_{2},\|\|\hat{\beta}\|_{2}^{-1} \hat{\beta}+\beta^{*}\|_{2}\right\}$ to measure the estimation error.\par
\begin{figure}[H]
\centering
\subfigure{
\begin{tikzpicture}
\pgfplotsset{every axis legend/.append style={at={(0.5,1)},
anchor=south}, every axis y label/.append style={at={(0.05,0.5)}}}
\begin{axis}[title={$(a)$: $f(x)=x^2$}, legend columns=1, legend pos=north east, title style={font=\tiny},xlabel = $\nu$,
    ylabel = {$\log(error)$},legend style={font=\tiny}, xtick={2.1,2.5,3,3.5,4}, font=\tiny,height=7.3cm, width=7.5cm, xmajorgrids=true]
\addplot[
    color=red,
    mark=triangle]
    coordinates{(2.1,-1.503259)(2.2,-1.559149)(2.3,-1.633823)(2.4,-1.644961)
(2.5,-1.669965)(2.6,-1.709761)(2.7,-1.732395)(2.8,-1.734232)(2.9,-1.768421)(3,-1.765554)
(3.1,-1.781875)(3.2,-1.816321)(3.3,-1.798378)(3.4,-1.824957)(3.5,-1.820158)(3.6,-1.817413)
(3.7,-1.838129)(3.8,-1.839817)(3.9,-1.847849)(4,-1.829449)};
\addplot[
    color=blue,
    mark=x]
    coordinates{(2.1,-1.63176)(2.2,-1.684010)(2.3,-1.690527)(2.4,-1.716955)
(2.5,-1.776271)(2.6,-1.765189)(2.7,-1.801093)(2.8,-1.819892)(2.9,-1.827590)(3,-1.835110)
(3.1,-1.869931)(3.2,-1.874254)(3.3,-1.889851)(3.4,-1.907318)(3.5,-1.918755)(3.6,-1.929415)
(3.7,-1.962345)(3.8,-1.943688)(3.9,-1.941375)(4,-1.954671)};
\legend{DP ,RDP}
\end{axis}
\end{tikzpicture}}
\subfigure{
\begin{tikzpicture}
\pgfplotsset{every axis legend/.append style={at={(0.5,1)},
anchor=south}, every axis y label/.append style={at={(0.05,0.5)}}}
\begin{axis}[title={$(b)$: $f(x)=x^4/2$}, legend columns=1,title style={font=\tiny},xlabel = $\nu$,
    ylabel = {$\log(error)$},legend style={at={(0.8,0.55)}, anchor=north,legend columns=-1}, xtick={2.1,2.5,3,3.5,4}, font=\tiny,height=7.3cm, width=7.5cm, xmajorgrids=true]
\addplot[
    color=red,
    mark=triangle]
    coordinates {
(2.1,-1.319440)(2.2,-1.340429)(2.3,-1.304078)(2.4,-1.328930)(2.5,-1.307862)(2.6,-1.319013)(2.7,-1.310395)
(2.8,-1.321727)(2.9,-1.331042)(3,-1.353881)(3.1,-1.335139)(3.2,-1.343897)(3.3,-1.344826)(3.4,-1.318884)
(3.5,-1.338638)(3.6,-1.327973)(3.7,-1.340681)(3.8,-1.351630)(3.9,-1.342414)(4,-1.337287)};
\addplot[
    color=blue,
    mark=x]
    coordinates {
(2.1,-1.678638)(2.2,-1.693967)(2.3,-1.714338)(2.4,-1.707085)(2.5,-1.739408)(2.6,-1.732043)(2.7,-1.763208)
(2.8,-1.750295)(2.9,-1.767084)(3,-1.787261)(3.1,-1.781135)(3.2,-1.81118)(3.3,-1.802482)(3.4,-1.808137)
(3.5,-1.81269)(3.6,-1.822672)(3.7,-1.824108)(3.8,-1.804472)(3.9,-1.827947)(4,-1.836986)};
\legend{DP ,RDP}
\end{axis}
\end{tikzpicture}}
\caption{Comparison between RDP and DP with varying $\nu$: $(a)$ the dimension $d=50$, the subsample size $n=500$ and $m=20$ with 100 replications $(b)$ the dimension $d=100$, the subsample size $n=800$ and $m=20$ with 50 replications.}
\label{fig5}
\end{figure}
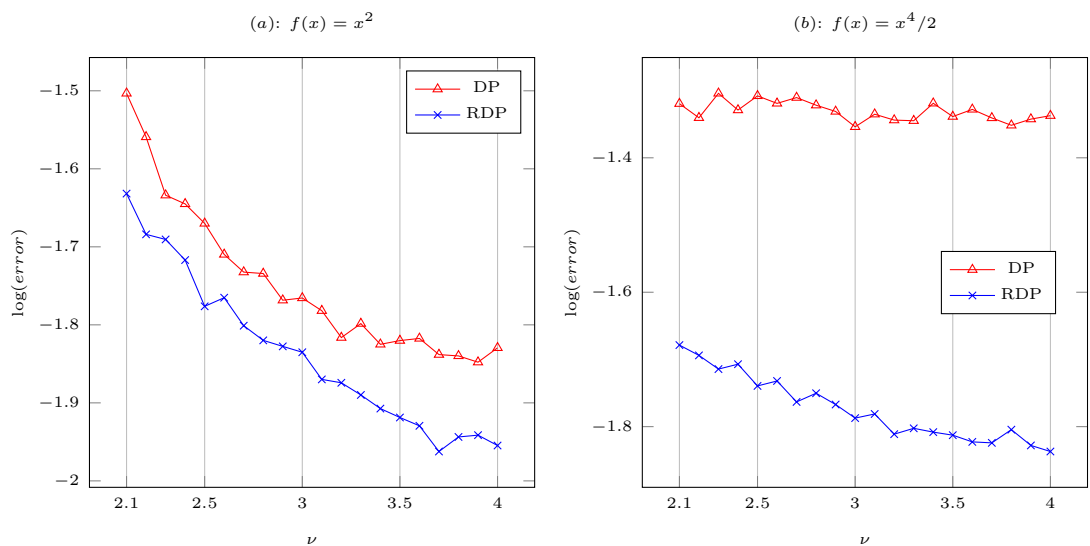
As shown in Figure \ref{fig5}, there is a more remarkable benefit in using RDP in the presence of heavy-tailed noise.\par
Finally, we choose the Swarm behaviour dataset\footnote{The original data consists of 24017 instances and 2400 features. We employ the first 1000 features of each instance as the sample.} (Dua and Graff (2017)\cite{UCI machine}) from the UCI machine learning repository as a real data example to test Algorithm \ref{Al}. We randomly select $N=24000$ instances and conduct PCA in a decentralized manner. The full sample covariance matrix is considered as the population covariance matrix. The performance of two algorithms is shown as below.
\begin{table}[H]
    \centering
    \fontsize{10}{10}\selectfont
\caption{The estimation error of two algorithms with different $K$ and $m$ when the experiments are repeated 50 times.}
\begin{tabular}{|c|c|c|c|c|c|}
\hline
\diagbox{$m$}{$K$}& 4 & 5 & 6 & 7 & \\
\cline{1-6} \multirow{2}{*}{$10$} & $-3.7060_{(0.034)}$ & $-2.6527_{(0.233)}$ & $-2.7707_{(0.080)}$ &$-1.2707_{(0.193)}$& RDP \\
\cline{2-6}  & $-3.7996_{(0.153)}$ & $-1.7987_{(0.200)}$ &$-1.7902_{(0.188)}$& $-0.4594_{(0.330)}$& DP \\
\cline{1-6} \multirow{2}{*}{$20$} & $-3.2003_{(0.029)}$ & $-1.9633_{(0.189)}$ & $-2.0868_{(0.096)}$ &$-0.6892_{(0.214)}$& RDP\\
\cline{2-6}  & $-2.4694_{(0.106)}$ & $-1.2009_{(0.189)}$ &$-1.1219_{(0.156)}$& $0.2297_{(0.124)}$&DP\\
\cline{1-6} \multirow{2}{*}{$30$} & $-2.8153_{(0.024)}$ & $-1.7639_{(0.133)}$ & $-1.6917_{(0.081)}$ &$-0.3158_{(0.232)}$& RDP\\
\cline{2-6}  & $-2.4487_{(0.111)}$ & $-0.7979_{(0.169)}$ &$-0.6484_{(0.202)}$& $0.3536_{(0.028)}$&DP
\\
\hline
\end{tabular}
\label{tab3}
\end{table}
Table \ref{tab3} shows that RDP affords a great improvement on the average estimation error.
\section{\textbf{Discussion}}\label{SecDiscussion}
In this paper, we extend the theoretical result established by Fan et al. (2019)\cite{fan2019distributed} to the case that bounded fourth moment of the sample with symmetric innovation and the bounded $6$-th moment for the asymmetric case. Numerical simulations indicate that the algorithm can adequately deal with heavy-tailed data and outliers. For the asymmetric heavy-tailed distribution, the experiment implies that the case of finite $4$-th moment still holds, but whether it can be further extended in theory requires further investigation. Moreover, the tightness of Theorem \ref{theorem1} is also unknown, which will be explored in future research.

\appendix
\section{\textbf{Appendix}}
\subsection{\textbf{Proof of Lemma \ref{lemma}}}
\begin{proof}
Instead of $\max(e^x-1,0)$ in the proof of Theorem 3.2 in Minsker (2018)\cite{Minsker_2018}, we define $\phi(x):= e^{x}-x-1$. For $t \geq 0,$
$$
\begin{array}{l}
\operatorname{P}\left(\lambda_{\max}\left(n\widehat{\Sigma}_{n}(\alpha,\theta)-n\Sigma\right) \geq nt\right)
=\operatorname{P}\left(\phi\left(\lambda_{\max} \left(n\theta\widehat{\Sigma}_{n}(\alpha,\theta)-n\theta\Sigma\right) \geq \phi(n\theta t)\right)\right) \\ \leq \frac{1}{\phi(n\theta t)}\mathbb{E} \operatorname{tr}\left[ \phi\left(n\theta\widehat{\Sigma}_{n}(\alpha,\theta)-n\theta\Sigma\right)\right].
\end{array}
$$
By following the proof of Lemma 3.1 in Minsker (2018)\cite{Minsker_2018}. we can obtain
$$
\mathbb{E} \operatorname{tr}\left[ \exp\left(n\theta\widehat{\Sigma}_{n}(\alpha,\theta)-n\theta\Sigma\right)\right]
\leq\operatorname{tr} \left[\exp \left(\sum_{i=1}^{n}c_{\alpha}\theta^{\alpha} \mathbb{E} |X_{i}X_{i}^{T}|^{\alpha} \right)\right].
$$
Due to
$$
\begin{array}{l}
-\log \left(I-\theta X_{i} {X_{i}}^{T}+c_{\alpha}\theta^{\alpha} {|X_{i} {X_{i}}^{T}|}^{\alpha}\right) \preceq \psi_{\alpha}\left(\theta X_{i} {X_{i}}^{T}\right)\\
\preceq \log \left(I+\theta X_{i} {X_{i}}^{T}+c_{\alpha}\theta^{\alpha} {|X_{i} {X_{i}}^{T}|}^{\alpha}\right)
\end{array}
$$
and $\log(1+x)\leq x$, it yields
$$\begin{array}{l}
-\mathbb{E}\operatorname{tr}\left[\left(n\theta\widehat{\Sigma}_{n}(\alpha,\theta)-n\theta\Sigma\right)\right]
\leq \mathbb{E}\operatorname{tr}\left[\sum_{i=1}^{n}\log \left(I-\theta X_{i} {X_{i}}^{T}+c_{\alpha}\theta^{\alpha} {|X_{i} {X_{i}}^{T}|}^{\alpha}\right)+n\theta\Sigma\right]\\
\leq \mathbb{E}\operatorname{tr}\left[\sum_{i=1}^{n}\left(-\theta X_{i} {X_{i}}^{T}+c_{\alpha}\theta^{\alpha} {|X_{i} {X_{i}}^{T}|}^{\alpha}\right)+n\theta\Sigma\right]=\operatorname{tr} \left[\sum_{i=1}^{n}c_{\alpha}\theta^{\alpha} \mathbb{E} |X_{i}X_{i}^{T}|^{\alpha}\right].
\end{array}
$$
Therefore,
$$\mathbb{E} \operatorname{tr}\left[ \phi\left(n\theta\widehat{\Sigma}_{n}(\alpha,\theta)-n\theta\Sigma\right)\right]\leq \operatorname{tr} \left[\exp \left(\sum_{i=1}^{n}c_{\alpha}\theta^{\alpha} \mathbb{E} |X_{i}X_{i}^{T}|^{\alpha} \right)+\sum_{i=1}^{n}c_{\alpha}\theta^{\alpha} \mathbb{E} |X_{i}X_{i}^{T}|^{\alpha}-I\right].$$
Because $\frac{e^{x}-x-1}{x}=\sum_{i=1}^{\infty} \frac{x^{i}}{(i +1)!}$, we have
$$
\begin{array}{l}
\operatorname{tr} \left[\exp \left(\sum_{i=1}^{n}c_{\alpha}\theta^{\alpha} \mathbb{E} |X_{i}X_{i}^{T}|^{\alpha} \right)+\sum_{i=1}^{n}c_{\alpha}\theta^{\alpha} \mathbb{E} |X_{i}X_{i}^{T}|^{\alpha}-I\right]\\=\operatorname{tr}\left[c_{\alpha}\theta^{\alpha} \sqrt{n\mathbb{E} |X_{i}X_{i}^{T}|^{\alpha}}\left(2I+\sum_{k=2}^{\infty}{\left(c_{\alpha}\theta^{\alpha} n\mathbb{E} |X_{i}X_{i}^{T}|^{\alpha}\right)^{k-1}}/{k !}\right) \sqrt{n\mathbb{E} |X_{i}X_{i}^{T}|^{\alpha}}\right] \\
\leq \operatorname{tr}\left[c_{\alpha}\theta^{\alpha}n\mathbb{E} |X_{i}X_{i}^{T}|^{\alpha}\left(2+\sum_{k=2}^{\infty}{\left(c_{\alpha}\theta^{\alpha}n\left\|\mathbb{E} |X_{i}X_{i}^{T}|^{\alpha}\right\|_2\right)^{k-1}}/{k !}\right)\right]\\
= \bar{d}_{\alpha}\left(\exp \left(c_{\alpha}\theta^{\alpha}n\left\|\mathbb{E}
|X_{i}X_{i}^{T}|^{\alpha}\right\|_{2}\right)+c_{\alpha}\theta^{\alpha}n\left\|\mathbb{E} |X_{i}X_{i}^{T}|^{\alpha}\right\|_{2}-1\right).
\end{array}
$$
By $\frac{e^{x}}{e^{x}-x-1} \leq 1+\frac{2}{x}+\frac{2}{x^{2}}$ for $x>0$, it yields
$$\begin{array}{l}
\left(\exp \left(c_{\alpha}n\theta^{\alpha}\left\|\mathbb{E}
|X_{i}X_{i}^{T}|^{\alpha}\right\|_{2}\right)
+c_{\alpha}n\theta^{\alpha}\left\|\mathbb{E}
|X_{i}X_{i}^{T}|^{\alpha}\right\|_{2}-1\right)/\phi(n\theta t)\\
\leq \left(\exp \left(c_{\alpha}n\theta^{\alpha}\left\|\mathbb{E}
|X_{i}X_{i}^{T}|^{\alpha}\right\|_{2}-\theta nt\right)+c_{\alpha}n\theta^{\alpha}\left\|\mathbb{E}
|X_{i}X_{i}^{T}|^{\alpha}\right\|_{2}e^{-\theta nt}\right) \frac{e^{n\theta t}}{e^{n\theta t}-n\theta t-1}\\
\leq \left(\exp \left(c_{\alpha}n\theta^{\alpha}\left\|\mathbb{E}
|X_{i}X_{i}^{T}|^{\alpha}\right\|_{2}-n\theta t\right)+c_{\alpha}n\theta^{\alpha}\left\|\mathbb{E}
|X_{i}X_{i}^{T}|^{\alpha}\right\|_{2}e^{-\theta nt}\right)\left(1+\frac{2}{n\theta t}+\frac{2}{(n\theta t)^{2}}\right).
\end{array}
$$
Therefore,
\[
\operatorname{P}\left(\lambda_{\max}\left(\widehat{\Sigma}_{n}(\alpha,\theta)-\Sigma\right) \geq t\right) \leq \bar{d}_{\alpha}e^{-\theta nt}\left(e^{c_{\alpha}n \theta^{\alpha}v^{\alpha}}+c_{\alpha}n \theta^{\alpha}v^{\alpha}\right)\left(1+\frac{2}{\theta nt}+\frac{2}{(\theta nt)^2}\right).
\]
By the same way, we have
$$\operatorname{P}\left(\lambda_{\min}\left(\widehat{\Sigma}_{n}(\alpha,\theta)-\Sigma\right) \leq -t\right) \leq \bar{d}_{\alpha}e^{-\theta nt}\left(e^{c_{\alpha}n \theta^{\alpha}v^{\alpha}}+c_{\alpha}n \theta^{\alpha}v^{\alpha}\right)\left(1+\frac{2}{\theta nt}+\frac{2}{(\theta nt)^2}\right).$$
\end{proof}

\subsection{\textbf{Proof of Theorem \ref{theorem1}}}
\begin{proof}
For $\forall k, s \in [d]$ and $t>0$, by $\sigma_{k,s}:=\Sigma_{(k,s)}$, we have
$$
\begin{array}{l}
\operatorname{P}\left(\widehat{\sigma}_{k,s}-\sigma_{k,s}\geq t\right)=\operatorname{P}\left(\frac{n}{\tau_{k,s}}\widehat{\sigma}_{k,s}\geq \frac{n}{\tau_{k,s}}\sigma_{k,s}+\frac{nt}{\tau_{k,s}}\right)\leq e^{-\frac{n}{\tau_{k,s}}(t+\sigma_{k,s})}\mathbb{E}\left[\exp\left(n\widehat{\sigma}_{k,s}/\tau_{k,s}\right)\right]\\
=e^{-\frac{n}{\tau_{k,s}}(t+\sigma_{k,s})}\mathbb{E}\left[\exp\left(\sum_{i=1}^{n}
\psi_{\tau_{k,s}}\left({X_{i}}_{(k)}{X_{i}}_{(s)}\right)/ \tau_{k,s}\right)\right]\\=e^{-\frac{n}{\tau_{k,s}}(t+\sigma_{k,s})}\mathbb{E}\left[\exp\left(\sum_{i=1}^{n}
\psi_{1}\left({X_{i}}_{(k)}{X_{i}}_{(s)}/ \tau_{k,s}\right)\right)\right]\\
\leq e^{-\frac{n}{\tau_{k,s}}(t+\sigma_{k,s})}\mathbb{E}\left[\prod_{i=1}^{n}\left(1+{X_{i}}_{(k)}{X_{i}}_{(s)}/ \tau_{k,s}+\left|{X_{i}}_{(k)}{X_{i}}_{(s)}/ \tau_{k,s}\right|^{\alpha}\right)\right]\\
=e^{-\frac{n}{\tau_{k,s}}(t+\sigma_{k,s})}\prod_{i=1}^{n}\left(1+\sigma_{k,s}/ \tau_{k,s}
+\mathbb{E}\left|{X_{i}}_{(k)}{X_{i}}_{(s)}/ \tau_{k,s}\right|^{\alpha}\right)\\
\leq e^{-\frac{n}{\tau_{k,s}}(t+\sigma_{k,s})}\prod_{i=1}^{n}\exp\left(\sigma_{k,s}/ \tau_{k,s}+\mathbb{E}\left|{X_{i}}_{(k)}{X_{i}}_{(s)}/ \tau_{k,s}\right|^{\alpha}\right)\\
=\exp\left(n\mathbb{E}\left|{X_{i}}_{(k)}{X_{i}}_{(s)}\right|^{\alpha}/ \tau_{k,s}^{\alpha}-nt/\tau_{k,s}\right).
\end{array}
$$
Setting $\tau_{k,s}=\left(\frac{2\mathbb{E}\left|{X_{i}}_{(k)}{X_{i}}_{(s)}\right|^{\alpha}}{t}\right)^{\frac{1}{\alpha-1}}$, it yields
$$
\operatorname{P}\left(\widehat{\sigma}_{k,s}-\sigma_{k,s}\geq t\right)\leq \exp \left(-n\left(\mathbb{E}\left|{X_{i}}_{(k)}{X_{i}}_{(s)}\right|^{\alpha}\right)^{\frac{-1}{\alpha-1}}
\left(\frac{t}{2}\right)^{\frac{\alpha}{\alpha-1}}\right).
$$
When $t=2\left(\mathbb{E}\left|{X_{i}}_{(k)}{X_{i}}_{(s)}\right|^{\alpha}\right)^{\frac{1}{\alpha}}\left(\frac{2\log d -\log \delta}{n}\right)^{\frac{\alpha-1}{\alpha}}$, we have
$$
\operatorname{P}\left(\widehat{\sigma}_{k,s}-\sigma_{k,s}\geq 2\sqrt[\alpha]{\mathbb{E}\left|{X_{i}}_{(k)}{X_{i}}_{(s)}\right|^{\alpha}}\left(\frac{2\log d -\log \delta}{n}\right)^{\frac{\alpha-1}{\alpha}}\right)\leq \frac{\delta}{d^{2}}.
$$
Therefore,
$$
\operatorname{P}\left(|\widehat{\sigma}_{k,s}-\sigma_{k,s}|\geq 2\sqrt[\alpha]{\mathbb{E}\left|{X_{i}}_{(k)}{X_{i}}_{(s)}\right|^{\alpha}}\left(\frac{2\log d -\log \delta}{n}\right)^{\frac{\alpha-1}{\alpha}}\right)\leq \frac{2\delta}{d^{2}}.
$$
By the union bound, it yields
$$
\operatorname{P}\left(\left\|\widehat{\Sigma}-\Sigma\right\|_{\max}\geq 2M\left(\frac{2\log d -\log \delta}{n}\right)^{\frac{\alpha-1}{\alpha}}\right)\leq (1+d^{-1})\delta.
$$
\end{proof}

\subsection{\textbf{Proof of Lemma \ref{lemma4}}}
\begin{proof}
Define $D_{j}:=I-2 e_{j}e_{j}^{T},$ for $\forall j \in[d]$. Suppose that $\widehat{\lambda} \in \mathbb{R}$ and $\widehat{v} \in \mathbb{S}^{d-1}$ are an eigenvalue and the correspondent eigenvector of $$\widehat{\Sigma}_{n}(\alpha, \tau)=\frac{1}{n} \sum_{i=1}^{n} \psi_{\tau}\left(\left\|X_{i}^{(\ell)}\right\|_{2}^{2}\right) \frac{X_{i}^{(\ell)} X_{i}^{(\ell)T}}{\left\|X_{i}^{(\ell)}\right\|_{2}^{2}}$$
such that $\widehat{\Sigma}_{n}(\alpha, \tau) \hat{v}=\widehat{\lambda} \hat{v}$.
Let $\Sigma^{(\ell)}=V^{(\ell)} \Lambda^{(\ell)} V^{T(\ell)}$ be the eigendecomposition of $\Sigma^{(\ell)}$. For ease of notation, we remove the superscript $\ell$, and define $Z_{i}=\Lambda^{-\frac{1}{2}} V^{T} X_{i}$ and $\widehat{S}=\frac{1}{n } \sum_{i=1}^{n} \psi_{\tau}\left(\left\|X_{i}\right\|_{2}^{2}\right) \frac{Z_{i} Z_{i}^{T}}{\left\|X_{i}\right\|_{2}^{2}} . $ It yields $\widehat{\Sigma}_{n}(\alpha, \tau)=V \Lambda^{\frac{1}{2}} \widehat{S} {\Lambda}^{\frac{1}{2}} {V}^{T}$.
We denote the matrix $\check{\Sigma}:={V} {\Lambda}^{\frac{1}{2}} {D}_{j} \widehat{{S}} {D}_{j} {\Lambda}^{\frac{1}{2}} {V}^{T}$. Because $\left\{X_{i}\right\}_{i=1}^{n}$ is symmetric innovation, we have ${Z}_{i}\stackrel{d}{=}D_{j}{Z}_{i}:={{Z}_{i}}^{*}$,  and
\begin{align*}
\check{\Sigma}= {V} {\Lambda}^{\frac{1}{2}}\left(\frac{1}{n } \sum_{i=1}^{n} \psi_{\tau}\left(\left\|X_{i}\right\|_{2}^{2}\right) \frac{{Z_{i}}^{*} {{Z_{i}}^{*}}^{T}}{\left\|X_{i}\right\|_{2}^{2}}\right) {\Lambda}^{\frac{1}{2}} {V}^{T}
= \frac{1}{n} \sum_{i=1}^{n} \psi_{\tau}\left(\left\|{V} {\Lambda}^{\frac{1}{2}}Z_{i}\right\|_{2}^{2}\right) \frac{{V} {\Lambda}^{\frac{1}{2}}{Z_{i}}^{*} {{Z_{i}}^{*}}^{T} {\Lambda}^{\frac{1}{2}} {V}^{T}}{\left\|{V} {\Lambda}^{\frac{1}{2}}Z_{i}\right\|_{2}^{2}} .
\end{align*}
Note that $\left\|{V} {\Lambda}^{\frac{1}{2}}Z_{i}\right\|_{2}^{2}=\left\|{V} {\Lambda}^{\frac{1}{2}}{Z_{i}}^{*}\right\|_{2}^{2}.$ Hence, we have
$$
\check{\Sigma}=\frac{1}{n} \sum_{i=1}^{n} \psi_{\tau}\left(\left\|{V} {\Lambda}^{\frac{1}{2}}{Z_{i}}^{*}\right\|_{2}^{2}\right) \frac{{V} {\Lambda}^{\frac{1}{2}}{Z_{i}}^{*} {{Z_{i}}^{*}}^{T} {\Lambda}^{\frac{1}{2}} {V}^{T}}{\left\|{V} {\Lambda}^{\frac{1}{2}}{Z_{i}}^{*}\right\|_{2}^{2}},
$$
\begin{align*}
\widehat{\Sigma}_{n}(\alpha,\tau)= \frac{1}{n} \sum_{i=1}^{n} \psi_{\tau}\left(\left\|X_{i}\right\|_{2}^{2}\right) \frac{X_{i} X_{i}^{T}}{\left\|X_{i}\right\|_{2}^{2}}
= \frac{1}{n} \sum_{i=1}^{n} \psi_{\tau}\left(\left\|{V} {\Lambda}^{\frac{1}{2}}Z_{i}\right\|_{2}^{2}\right) \frac{{V} {\Lambda}^{\frac{1}{2}}Z_{i} {Z_{i}}^{T} {\Lambda}^{\frac{1}{2}} {V}^{T}}{\left\|{V} {\Lambda}^{\frac{1}{2}}Z_{i}\right\|_{2}^{2}}.
\end{align*}
Therefore, we get that $\widehat{\Sigma}_{n}(\alpha,\tau)$ and $\check{\Sigma}$ are identically distributed. The rest of the proof is the same as that of Theorem 2 in Fan et al. (2019)\cite{fan2019distributed}.
\end{proof}

\subsection{\textbf{Proof of Theorem \ref{theorem5}}}
\begin{proof}
By Lemma \ref{lemma1.5} and $x<e^{x}$, we can get that for $\tau=O\left(\sigma \cdot \sqrt{n}\right)$,
$$\operatorname{P}\left(\left\|\widehat{\Sigma}_{n}(2,\tau)-\Sigma\right\|_{2} \geq t\right) \leq C_{1} \bar{d}\left(1+\frac{2\sigma}{t\sqrt{n}}
+\frac{2\sigma^{2}}{t^{2}n}\right) \exp \left(-\frac{t\sqrt{n}}{\sigma}\right). $$
By the equivalent definition of sub-exponential random variable and $\psi_{1}$-norm,
$$
\begin{array}{l}
\mathbb{E}\left\|\widehat{\Sigma}_{n}(2,\tau)-\Sigma\right\|_{2}^{k}
=  \int_{0}^{\infty} \operatorname{P}\left(\left\|\widehat{\Sigma}_{n}(2,\tau)-\Sigma\right\|_{2}^{k}>s\right) \mathrm{d} s \\
= \int_{0}^{\infty} \operatorname{P}\left(\left\|\widehat{\Sigma}_{n}(2,\tau)-\Sigma\right\|_{2}>s^{1 / k}\right) \mathrm{d} s \\
\leq \int_{0}^{\infty} C_{1} \bar{d}\left(1+\frac{2\sigma}{s^{1/k}\sqrt{n}}
+\frac{2\sigma^{2}}{s^{2/k}n}\right) \exp \left(-\frac{s^{1/k}\sqrt{n}}{\sigma}\right)\mathrm{d}s
\\
= C_{1}\bar{d} \left(\frac{\sigma}{\sqrt{n}}\right)^{k} k \int_{0}^{\infty} e^{-u}\left( u^{k-1}+ 2 u^{k-2} + 2u^{k-3}\right) \mathrm{d} u  \;\; \;\; \left(u:=\frac{s^{1/k}\sqrt{n}}{\sigma}\right)\\
= C_{1}\bar{d}\left(\frac{\sigma}{\sqrt{n}}\right)^{k} k \left(\Gamma(k)+2\Gamma(k-1)+2\Gamma(k-2)\right).
\end{array}
$$
Because $\Gamma(k) \leq k^{k}$ and for any $k \geq 1$, $k^{1 / k} \leq e^{1 / e} \leq 2$, we have
\[
\left(C_{1}\bar{d}\left(\frac{\sigma}{\sqrt{n}}\right)^{k} k \left(\Gamma(k)+2\Gamma(k-1)+2\Gamma(k-2)\right)\right)^{1 / k} \leq C\left(\bar{d}\right)^{1/k}{\frac{\sigma}{\sqrt{n}}} k
\leq C\bar{d}{\frac{\sigma}{\sqrt{n}}} k.
\]
\\
Hence,
$\left\|\left\|\widehat{\Sigma}_{n}(2,\tau)-\Sigma\right\|_{2}\right\|_{\psi_{1}}=\sup _{k \geq 1}\left(\mathbb{E}\left\|\widehat{\Sigma}_{n}(2,\tau)-\Sigma\right\|_{2}^{k}\right)^{1/k}/k\leq C\bar{d}{\frac{\sigma}{\sqrt{n}}}.$\\
By the Davis-Kahan theorem\cite{10.1093/biomet/asv008},
\begin{equation}\label{eq6}
\left\|\rho\left(\widetilde{{V}}_{K},  {V}_{K}\right)\right\|_{\psi_{1}}\lesssim \left\|\left\|\widetilde{{\Sigma}}-{V}_{K} {V}_{K}^{T}\right\|_{F}\right\|_{\psi_{1}} \leq \left\|\left\|\widetilde{{\Sigma}}-{\Sigma}^{*}\right\|_{F}\right\|_{\psi_{1}}+\left\|{\Sigma}^{*}-{V}_{K} {V}_{K}^{T}\right\|_{F}.
\end{equation}
By the robust covariance version of Lemma 1 and Theorem 2 in Fan et al. (2019)\cite{fan2019distributed}, if for all $\ell \in [m]$, $\|\mathbb{E}\widehat{{V}}_{K}^{(\ell)} \widehat{{V}}_{K}^{(\ell) T}-{V}_{K} {V}_{K}^{T}\|_{2} \leq 1 / 4$, the first term in (\ref{eq6}) can be written as
\begin{equation}\label{eq7}
\begin{aligned}
&\left\|\left\|\widetilde{{\Sigma}}-{\Sigma}^{*}\right\|_{F}\right\|_{\psi_{1}}
\lesssim \frac{1}{m} \sqrt{\sum_{\ell=1}^{m}\left\|\left\|\widehat{{V}}_{K}^{(\ell)} \widehat{{V}}_{K}^{(\ell) T}-\mathbb{E}\widehat{{V}}_{K}^{(\ell)} \widehat{{V}}_{K}^{(\ell) T}\right\|_{F}\right\|_{\psi_{1}}^{2}} \\&
\lesssim \frac{1}{m} \sqrt{\sum_{\ell=1}^{m}\left(\bar{d}_{(\ell)}\frac{\sqrt{K}}{\Delta_{(\ell)}} {\frac{\sigma_{(\ell)}}{\sqrt{n}}}\right)^{2}}\lesssim  \sqrt{\frac{1}{m}\sum_{\ell=1}^{m}\left(\bar{d}_{(\ell)} \frac{\sigma_{(\ell)}}{\Delta_{(\ell)}}\right)^{2}} \sqrt{\frac{K}{N}}.
\end{aligned}
\end{equation}
Since
\begin{equation}\label{eq8}
\begin{aligned}
&\left\|\mathbb{E}\widehat{{V}}_{K}^{(\ell)} \widehat{{V}}_{K}^{(\ell) T}-V_{K} V_{K}^{T}\right\|_{2}\leq \mathbb{E}\left\|\widehat{{V}}_{K}^{(\ell)} \widehat{{V}}_{K}^{(\ell) T}-V_{K} V_{K}^{T}\right\|_{2}\leq \mathbb{E}\left\|\widehat{{V}}_{K}^{(\ell)} \widehat{{V}}_{K}^{(\ell) T}-V_{K} V_{K}^{T}\right\|_{F}\\ &\leq \frac{\sqrt{K}}{\Delta_{(\ell)}}\mathbb{E}\left\|\widehat{\Sigma}^{(\ell)}(2,\tau_{(\ell)})-\Sigma^{(\ell)}\right\|_{2} \leq
\frac{\sqrt{K}}{\Delta_{(\ell)}}\left\|\left\|\widehat{\Sigma}^{(\ell)}(2,\tau_{(\ell)})-\Sigma^{(\ell)}\right\|_{2}\right\|_{\psi_{1}}
\lesssim \bar{d}_{(\ell)}\frac{\sqrt{K}}{\Delta_{(\ell)}}\frac{\sigma_{(\ell)}}{\sqrt{n}},
\end{aligned}
\end{equation}
we obtain that if $C_{1}$ is sufficiently large such that $n \geq C_{1}K \max_{\ell \in [m]}\left(\bar{d}_{(\ell)}\frac{\sigma_{(\ell)}}{\Delta_{(\ell)}}\right)^2$, (\ref{eq8}) implies that $\| \mathbb{E}\widehat{{V}}_{K}^{(\ell)} \widehat{{V}}_{K}^{(\ell) T}-{V}_{K} {V}_{K}^{T}\|_{2} \leq 1 / 4< 1/2$ for all $\ell \in [m]$. Therefore, by (\ref{eq7}) and Lemma \ref{lemma4}, we have for some constant $C_{2}$,
$$\left\|\rho\left(\widetilde{{V}}_{K},  {V}_{K}\right)\right\|_{\psi_{1}} \leq C_{2}\sqrt{\frac{1}{m}\sum_{\ell=1}^{m}\left(\bar{d}_{(\ell)} \frac{\sigma_{(\ell)}}{\Delta_{(\ell)}}\right)^{2}} \sqrt{\frac{K}{N}}.$$
\end{proof}

\subsection{\textbf{Proof of Lemma \ref{lemma4.5}}}
\begin{proof}
For $\forall v \in \mathcal{S}^{d-1}$, we have
$$
\begin{array}{l}
\left\|\mathbb{E}|X_{i} X_{i}^{T}|^{2}\right\|_{2}=\mathbb{E}\left(\left\|X_{i}\right\|_{2}^{2}(v^{T}X_{i})^2\right)
=\sum_{j=1}^{d}\mathbb{E}\left(x_{ij}^{2}(v^{T}X_{i})^2\right)\\ \leq \sum_{j=1}^{d}\left(\mathbb{E}|x_{ij}|^{3}\right)^{\frac{2}{3}}\left(\mathbb{E}(v^{T}X_{i})^6\right)^{\frac{1}{3}}
\leq\sum_{j=1}^{d}\left(\mathbb{E}x_{ij}^{6}\right)^{\frac{1}{3}}\left(\mathbb{E}(v^{T}X_{i})^6\right)^{\frac{1}{3}}
\leq d R_{(1)}^{\prime\frac{2}{3}}<\infty.
\end{array}
$$
Therefore, define ${\Omega}=\widehat{\Sigma}_{n}^{(1)}(2, \tau_{(1)})-{\Sigma}^{(1)}$,  ${\Gamma}={V}_{K} {V}_{K}^{T}$,  $\widehat{{\Gamma}}=\widehat{{V}}_{K}^{(1)} \widehat{{V}}_{K}^{(1) T}$, ${\Theta}=f\left({\Omega V}_{K}\right) {V}_{K}^{T}+{V}_{K} f\left({\Omega V}_{K}\right)^{T}$ where $f$ is a linear function defined in Lemma 2 of Fan et al. (2019)\cite{fan2019distributed}, ${\Phi} =\widehat{{\Gamma}}-{\Gamma}-{\Theta}$ and $\omega=\|{\Omega}\|_{2} / \Delta$. Since
$$
\widehat{{\Gamma}}-{\Gamma}={\Phi} 1_{\{\omega \leq 1 / 10\}}+(\widehat{{\Gamma}}-{\Gamma}) 1_{\{\omega>1 / 10\}}-{\Theta} 1_{\{\omega>1 / 10\}}+\Theta,
$$
we have
\begin{equation}\label{eq3}
\|\mathbb{E} \widehat{{\Gamma}}-{\Gamma}\|_{F} \leq \mathbb{E}\left(\|{\Phi}\|_{F} 1_{\{\omega \leq 1 / 10\}}\right)+\mathbb{E}\left(\|\widehat{{\Gamma}}-{\Gamma}\|_{F} 1_{\{\omega>1 / 10\}}\right)+\mathbb{E}\left(\|{\Theta}\|_{F} 1_{\{\omega > 1 / 10\}}\right)+\|\mathbb{E}\left(\Theta\right)\|_{F}.
\end{equation}
By Theorem 3 in Fan et al. (2019)\cite{fan2019distributed}, it yields
\begin{equation}\label{eq4}
\mathbb{E}\left(\|{\Phi}\|_{F} 1_{\{\omega \leq 1 / 10\}}\right)+\mathbb{E}\left(\|\widehat{{\Gamma}}-{\Gamma}\|_{F} 1_{\{\omega>1 / 10\}}\right)+\mathbb{E}\left(\|{\Theta}\|_{F} 1_{\{\omega > 1 / 10\}}\right) \lesssim \sqrt{K} \mathbb{E}\omega^2.
\end{equation}
Since $\mathbb{E}\left(\Omega\right)=\mathbb{E}\left(\psi_{\tau_{(1)}}\left(\left\|X_{i}\right\|_{2}^{2}\right) \frac{X_{i} X_{i}^{T}}{\left\|X_{i}\right\|_{2}^{2}}-X_{i} X_{i}^{T}\right)=\mathbb{E}\left(\left(\psi_{\tau_{(1)}}\left(\left\|X_{i}\right\|_{2}^{2}\right) /\left\|X_{i}\right\|_{2}^{2}-1\right)X_{i} X_{i}^{T}\right)$, for $\forall v \in \mathcal{S}^{d-1}$, we have
\begin{align*}
\mathbb{E}&\left(\left(\psi_{\tau_{(1)}}\left(\left\|X_{i}\right\|_{2}^{2}\right) /\left\|X_{i}\right\|_{2}^{2}-1\right)v^{T}X_{i} X_{i}^{T}v\right)=\mathbb{E}\left(\frac{\tau_{(1)}-\left\|X_{i}\right\|_{2}^{2}} {\left\|X_{i}\right\|_{2}^{2}}(v^{T}X_{i})^2 1_{\{\left\|X_{i}\right\|_{2}^{2}>\tau_{(1)}\}}\right)\\
\leq & \mathbb{E}\left((v^{T}X_{i})^2 1_{\{\left\|X_{i}\right\|_{2}^{2}>\tau_{(1)}\}}\right)\leq
\left(\mathbb{E}(v^{T}X_{i})^6\right)^{\frac{1}{3}}\left(\mathrm{P}\left(\left\|X_{i}
\right\|_{2}^{2}>\tau_{(1)}\right)\right)^{\frac{2}{3}}\\ \leq & R_{(1)}^{\prime\frac{1}{3}}\left(\mathbb{E}\left\|X_{i}\right\|_{2}^{6}/\tau_{(1)}^3\right)^{\frac{2}{3}}
\lesssim R_{(1)}^{\prime}d^{2}/(\sigma_{(1)}^2 n)
\end{align*}
where the second and third inequalities follow from H{\"o}lder and Markov inequality. The last inequality follows from $C_{r}$ inequality. Hence, $\|\mathbb{E}\left(\Omega\right)\|_{2}\lesssim R_{(1)}^{\prime}d^{2}/\left(\sigma_{(1)}^2 n\right)$ and
\begin{equation}\label{eq10}
\|\mathbb{E}\left(\Theta\right)\|_{F}\lesssim \left\|f\left({\mathbb{E}\left(\Omega\right) V}_{K}\right)\right\|_{F}\leq \sqrt{K}\|\mathbb{E}\left(\Omega\right)\|_{2} / \Delta_{(1)} \lesssim \frac{\sqrt{K}}{\Delta_{(1)}}\frac{R_{(1)}^{\prime}d^{2}}{\sigma_{(1)}^2 n}.
\end{equation}
Finally, combing (\ref{eq3})-(\ref{eq10}), it can be shown that
\begin{align*}
\|\mathbb{E} \widehat{{\Gamma}}-{\Gamma}\|_{F} \lesssim & \sqrt{K} \mathbb{E} \omega^{2}+\|\mathbb{E}\left(\Theta\right)\|_{F}\lesssim  \sqrt{K} \Delta^{-2}\|\| {\Omega}\left\|_{2}\right\|_{\psi_{1}}^{2}+\frac{\sqrt{K}}{\Delta_{(1)}}\frac{R_{(1)}^{\prime}d^{2}}{\sigma_{(1)}^2 n} \\ \lesssim & \left(\bar{d}_{(1)}\frac{\sigma_{(1)}}{\Delta_{(1)}}\right)^2\frac{\sqrt{K}}{n}+
\frac{R_{(1)}^{\prime}d^{2}}{\sigma_{(1)}^2\Delta_{(1)}}\frac{\sqrt{K}}{n}.
\end{align*}
\end{proof}

\subsection{\textbf{Proof of Theorem \ref{theorem6}}}
\begin{proof}
By Lemma \ref{lemma4.5} and (\ref{eq7}), we obtain that when $n \geq C_{1}K \max_{\ell \in [m]}\left(\bar{d}_{(\ell)}\frac{\sigma_{(\ell)}}{\Delta_{(\ell)}}\right)^2$,
$$
\begin{aligned}
&\left\|\rho\left(\widetilde{{V}}_{K},  {V}_{K}\right)\right\|_{\psi_{1}}\leq \left\|\rho\left(\widetilde{{V}}_{K},  {V}_{K}^{\star}\right)\right\|_{\psi_{1}}+\rho\left({V}_{K}^{\star},  {V}_{K}\right)\\ &\lesssim  \sqrt{\frac{1}{m}\sum_{\ell=1}^{m}\left(\bar{d}_{(\ell)} \frac{\sigma_{(\ell)}}{\Delta_{(\ell)}}\right)^{2}} \sqrt{\frac{K}{N}}+ \frac{1}{m}\sum_{\ell=1}^{m}\left(\left(\bar{d}_{(\ell)}\frac{\sigma_{(\ell)}}{\Delta_{(\ell)}}\right)^{2}
+\frac{R_{(\ell)}^{\prime}d^{2}}{\sigma_{(\ell)}^2\Delta_{(\ell)}}\right)\frac{\sqrt{K}}{n}.
\end{aligned}
$$
Therefore, when the requirement on $m$ and $n$ is satisfied, we have for constants $C_{3}$,
$$
\left\|\rho\left(\widetilde{{V}}_{K},  {V}_{K}\right)\right\|_{\psi_{1}}\leq C_{3}\sqrt{\frac{1}{m}\sum_{\ell=1}^{m}\left(\bar{d}_{(\ell)} \frac{\sigma_{(\ell)}}{\Delta_{(\ell)}}\right)^{2}} \sqrt{\frac{K}{N}}.
$$
\end{proof}
%----------------------------------------------------------------------------------------
%	REFERENCE LIST
%----------------------------------------------------------------------------------------

\end{document}